\documentclass[12pt]{amsart}

\usepackage[
ps2pdf,                
colorlinks=true,
linkcolor=blue,
citecolor=green,
urlcolor=magenta
]{hyperref}

\usepackage{epsfig}
\usepackage{latexsym}
\usepackage{psfrag}
\usepackage{amssymb}
\usepackage{amsmath}
\usepackage{amsthm} 

\newcommand{\etc}{{\it etc}}
\newcommand{\ie}{{\it i.e.}}
\newcommand{\eg}{{\it e.g.}}

\renewcommand{\th}{{\rm th}}

\newcommand{\calF}{{\mathcal{F}}}
\newcommand{\calH}{{\mathcal{H}}}
\newcommand{\calM}{{\mathcal{M}}}

\newcommand{\RR}{{\mathbb{R}}}

\newcommand{\ZZ}{{\mathbb{Z}}}
\newcommand{\NN}{{\mathbb{N}}}
\newcommand{\BB}{{\mathbb{B}}}
\newcommand{\DD}{{\mathbb{D}}}

\renewcommand{\setminus}{{\smallsetminus}}

\newcommand{\image}{{\operatorname{image}}}

\newcommand{\bigmid}{{~\mid~}}
\newcommand{\st}{{\bigmid}}
\newcommand{\from}{{\colon}} 

\newcommand{\homeo}{{\medspace \cong \medspace}} 
 
\newcommand{\cross}{{\times}}

\newcommand{\closure}[1]{{\overline{#1}}}

\newcommand{\euler}{{\operatorname{\chi}}} 
\newcommand{\neigh}{{\operatorname{\eta}}} 
\newcommand{\bdy}{{\partial}} 

\newcommand{\interior}{{\operatorname{interior}}}
\newcommand{\connect}{{\#}} 

\newcommand{\RRPP}{{\mathbb{RP}}} 

\newcommand{\vbdy}{{\partial_v}} 
\newcommand{\hbdy}{{\partial_h}} 

\theoremstyle{plain}
\newtheorem{theorem}{Theorem}[section]
\newtheorem{corollary}[theorem]{Corollary}
\newtheorem{lemma}[theorem]{Lemma}

\theoremstyle{definition}
\newtheorem*{define}{Definition}

\newtheorem*{proofclaim}{Claim}

\newtheorem{remark}[theorem]{Remark}

\newtheorem*{question}{Question}

\newsavebox{\savepar}

\newcommand{{\NP}}{{\bf NP}}
\newcommand{{\size}}{{\rm size}} 
\newcommand{{\weight}}{{\rm weight}} 

\newcommand{\normal}[1]{{\widetilde{#1}}} 
\newcommand{\blocky}[1]{{\widehat{#1}}} 

\newcommand{{\core}}{{\rm core}} 
\newcommand{{\product}}{{\rm product}} 

\begin{document}

\title{Sphere recognition lies in NP}
\author{Saul Schleimer}
\address{\hskip-\parindent
        Saul Schleimer\\
        Department of Mathematics, UIC\\
        851 South Morgan Street\\
        Chicago, Illinois 60607}
\email{saul@math.uic.edu}

\date{\today}

\begin{abstract}
We prove that the three-sphere recognition problem lies in the
complexity class NP.  Our work relies on Thompson's original proof
that the problem is decidable [Math. Res. Let., 1994], Casson's
version of her algorithm, and recent results of Agol, Hass, and
Thurston [ArXiv, 2002].
\end{abstract}
\maketitle

\section{Introduction}

The three-sphere recognition problem asks: given a triangulation $T$,
is the underlying space $|T|$ homeomorphic to the three-sphere?
Clearly there are trivial points to check, such as whether $|T|$ is a
three-manifold and if $|T|$ is closed, connected, orientable, \etc.
After dealing with these issues:

\vspace{2mm}
\noindent
{\bf 
Theorem~\ref{Thm:CassonAlg}(Thompson~\cite{Thompson94},
Casson~\cite{Casson97}).
}
{\em 
There is an exponential time algorithm which, given a triangulation $T$,
decides whether or not $|T|$ is homeomorphic to the three-sphere.
}
\vspace{1mm}

Our goal is to show:

\vspace{2mm}
\noindent
{\bf Theorem~\ref{Thm:RecognitionInNP}.}
{\em 
The three-sphere recognition problem lies in the complexity class \NP.
}
\vspace{1mm}

That is, if $T$ is a triangulation of the three-sphere then there is a
polynomial sized proof of this fact.  The essential details of such a
proof is called a {\em certificate}.

As a corollary of Theorem~\ref{Thm:RecognitionInNP}:

\begin{corollary}
\label{Cor:ThreeBallInNP}
The three-ball recognition problem lies in the complexity class \NP.
\end{corollary}

To prove Corollary~\ref{Cor:ThreeBallInNP} we must produce a
certificate for every triangulation $T$ of the three-ball.  (See
Section~\ref{Sec:Prelim} for how to verify, in polynomial time, that
$T$ is a three-manifold.)  First write down a polynomial sized proof
that $S = \bdy T$ is two-sphere.  Next we give the certificate that
$D(T)$, the {\em double} of $T$, is a three-sphere using
Theorem~\ref{Thm:RecognitionInNP}.  This completes the certificate as,
by Alexander's Theorem (see Theorem~1.1 of~\cite{Hatcher01}), the
two-sphere $S$ cuts $D(T)$ into a pair of three-balls.


Theorem~\ref{Thm:RecognitionInNP} could be used to place other
problems in \NP: for example, the problem of deciding whether or not a
triangulated three-manifold is the complement of a knot in $S^3$.
Jaco and Sedgwick~\cite{JacoSedgwick03} show that in any nontrivial
knot complement
there is a meridional boundary slope which is not too long.  Combining
this with our work and a discussion on layered solid tori gives a
proof.


It is possible that these techniques also show that the lens space
recognition and surface bundle recognition problems lie in \NP.  (See
Section~\ref{Sec:Questions}.)  Corollary~\ref{Cor:ThreeBallInNP} may also
find application for certifying that a triangulated three-manifold is
Haken.


Before abandoning the introduction we mention that many other problems
in three-manifold topology lie in \NP.  Hass, Lagarias, and
Pippenger~\cite{HassEtAl99} have shown that the unknotting problem,
first solved by Haken, lies in \NP.  Agol, Hass and
Thurston~\cite{AgolEtAl02} have shown that the 3-manifold knot genus
problem is in fact \NP-complete.  (The 3-manifold knot genus problem
asks; given a genus $g$ and a knot $K$ in the one-skeleton of a
triangulated three-manifold $M$, is there an orientable spanning
surface for $K$ with genus at most $g$?) Also,
S.~Ivanov~\cite{Ivanov01} has shown that the three-sphere recognition
algorithm is in \NP~when the triangulations considered are assumed to
be {\em zero-efficient}.  (See Remark~\ref{Rem:VertexFundamental}.)

In this paper I rely heavily on material discussed
in~\cite{HassEtAl99} and~\cite{AgolEtAl02} as well as work of Jaco and
Rubinstein~\cite{JacoRubinstein02a} and a talk of Casson's at
MSRI~\cite{Casson97}.  I thank Ian Agol for suggesting that I tackle
this problem and for suggesting that ideas from my
thesis~\cite{Schleimer01} might be useful in its solution.

\section{Sketch of the main theorem}

Before diving into the details of the definitions we give a sketch of
Theorem~\ref{Thm:RecognitionInNP}.  We closely follow Casson's
algorithm for recognizing the three-sphere.

Fix $T$, a triangulation of $S^3$.  Produce a certificate $\{ (T_i,
v(S_i)) \}_{i = 0}^n$ as follows: The triangulation $T_0$ is equal to
$T$.  For every $i$ use Lemma~\ref{Lem:FundamentalSphereAlg} to find
$S_i$, a normal two-sphere in $T_i$ which is not vertex linking, if
such exists.  If $T$ is zero-efficient then
Lemma~\ref{Lem:FundamentalSphereAlg} provides $S_i$, an almost normal
two-sphere in $T_i$.  (Here $v(S_i)$ is a {\em surface vector}; a
concise representation of $S_i$.  See Section~\ref{Sec:NormalSurfaces}
for the definition.)

If $S_i$ is normal apply Theorem~\ref{Thm:CrushingAlg}: $T_{i+1}$ is
obtained from $T_i$ by crushing $T_i$ along $S_i$.  Briefly, we cut
$|T_i|$ along $S_i$, cone the resulting two-sphere boundary components
to points, and collapse non-tetrahedral cells of the resulting cell
structure to obtain the triangulation $T_{i+1}$.  This is discussed in
Section~\ref{Sec:Crushing}, below.

If $S_i$ is almost normal then obtain $T_{i+1}$ from $T_i$ by deleting
the component of $|T_i|$ which contains $S_i$.  Finally, the last
triangulation $T_n$ is empty, as is $S_n$.

That completes the construction of the certificate.  We now turn to
the procedure for checking a given certificate: that is, we cite a
series of polynomial time algorithms which verify each part of the
certificate.  First, check if $T = T_0$ using
Theorem~\ref{Thm:IdenticalAlg}.  Next, verify that $S_i$ is in fact
the desired surface by checking its Euler characteristic (see
Lemma~\ref{Lem:EulerAlg}) and checking that it is connected (see
Theorem~\ref{Thm:ConnectedAlg}).  Next, if $S_i$ is normal verify that
the triangulation $T_{i+1}$ is identical to the triangulation obtained
by crushing $T_i$ along $S_i$.  This uses
Theorem~\ref{Thm:CrushingAlg} and again uses
Theorem~\ref{Thm:IdenticalAlg}.  If $S_i$ is almost normal then check
that the component $T'$ of $T_i$ containing $S_i$ has $|T'| \homeo
S^3$ using Theorem~\ref{Thm:NewNormalizationAlg} and
Theorem~\ref{Thm:NormalizingS2}.

Finally, by Theorem~\ref{Thm:CrushingAndConnectSum}, for every $i$ we
have that $\connect |T_i| \homeo \connect |T_{i+1}|$ where the connect
sum on the left hand side ranges over the components of $|T_{i}|$
while the right hand side ranges over the components of $|T_{i+1}|$.
By definition the empty connect sum is $S^3$, and this finishes the
verification of the certificate.

\begin{remark}
\label{Rem:VertexFundamental}
As an aside, note that there is a special type of fundamental normal
(or almost normal) surface called a {\em vertex fundamental} surface.
These lie on extremal rays for a certain linear cone of embedded
normal (almost normal) surfaces.  They are thus vertices of the
projectivization of the cone.

It is possible to certify that such surfaces are, in fact, vertex
surfaces.  This in turn implies that they are fundamental.  Hass,
Lagarias, and Pippenger~\cite{HassEtAl99} use this fact to certify
connectedness of a normal surface and thus to show that unknotting is
in \NP.
This is markedly different from general fundamental surfaces
where it is not currently known how to certify that they are
fundamental.

S.~Ivanov has raised the possibility of using vertex fundamental
two-spheres in order to certify that the $S_i$ are connected.  We have
preferred to use Theorem~\ref{Thm:ConnectedAlg} which relies
on~\cite{AgolEtAl02}.  This is because we use~\cite{AgolEtAl02} in an
even more substantial way in Theorem~\ref{Thm:NormalizationAlg}.  I
would not need Theorem~\ref{Thm:NormalizationAlg} if I knew how to
certify that a triangulation is zero-efficient.  
\end{remark}

\section{Definitions}
\label{Sec:Define}

We begin with a naive discussion of complexity theory.  Please
consult~\cite{GareyJohnson79} for a more through treatment.


A {\em problem} $P$ is a function from a set of finite binary strings,
the {\em instances}, to another set of finite binary strings, the {\em
answers}.  A problem $P$ is a {\em decision problem} if the range of
$P$ is the set $\{0, 1\}$.  The length of a binary string in the
domain of $P$ is the {\em size} of the instance.  A {\em solution} for
$P$ is a Turing machine $\calM$ which, given input $T$ in its tape,
computes $P(T)$ and halts with only that output on its tape.  We will
engage in the usual abuse of calling such a Turing machine an {\em
algorithm} or a {\em procedure}.

An algorithm is {\em polynomial time} if there is a polynomial $q
\from \RR \to \RR$ so that the Turing machine $\calM$ finishes
computing in time at most $q(\size(T))$.  Computing bounds on $q$, or
even its degree, is often a difficult problem.  We will follow
previous treatments in algorithmic topology and leave this problem
aside.

A decision problem $P$ lies in the complexity class \NP~if there is a
polynomial $q' \from \RR \to \RR$ with the following property: For all
instances $T$ with $P(T) = 1$ there is proof of length at most
$q'(\size(T))$ that $P(T) = 1$.  Such a polynomial length proof is a
{\em certificate} for $T$.  More concretely: Suppose that there is a
polynomial $q''$ and a Turing machine $\calM''$ so that, for every
instance $T$ with $P(T) = 1$, there is a string $C$ where $\calM''$
run on $(T, C)$ outputs the desired proof that $P(T) = 1$ in time less
than $q''(\size(T))$.  Then, again, the problem $P$ is in \NP~and we
also call $C$ a certificate for $T$.

\begin{figure}
\psfrag{0}{$0$}
\psfrag{1}{$1$}
\psfrag{2}{$2$}
\psfrag{3}{$3$}
$$\begin{array}{c}
\epsfig{file=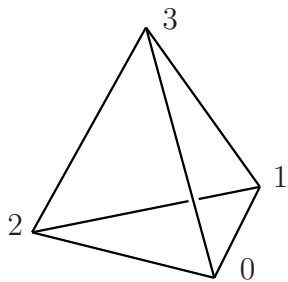, height = 3.5 cm}
\end{array}$$
\caption{A regular Euclidean tetrahedron with all side-lengths equal
  to one.}
\label{Fig:ModelTetrahedron}
\end{figure}

We now turn to topological considerations.  A {\em model tetrahedron}
$\tau$ is a copy of the regular Euclidean tetrahedron of side length
$1$ with vertices labelled by $0$, $1$, $2$, and $3$.  See
Figure~\ref{Fig:ModelTetrahedron} for a picture.  Label the six
edges by their vertices $(0,1)$, $(0,2)$, \etc.  Label the four
faces by the number of the vertex they do {\em not} contain.  The
standard orientation on $\RR^3$ induces an orientation on the model
tetrahedron which in turn induces orientations on the faces.

A {\em labelled triangulation} of size $n$ is a collection of $n$
model tetrahedra $\{\tau_i\}_{i = 1}^n$, each with a unique name, and
a collection of {\em face pairings}.  Here a face pairing is a triple
$( i, j, \sigma)$ specifying a pair of tetrahedra $\tau_i$ and
$\tau_j$ as well as an orientation reversing isometry $\sigma$ from a
face of $\tau_i$ to a face of $\tau_j$.  We will omit the labellings
when they are clear from context.

A triangulation is not required to be a simplical complex.  However
every face must appear in exactly two face pairings or in none.  Also,
no face may be glued to itself.  We do not require for a face pairing
$(i, j, \sigma)$ that $i \neq j$.


Let $|T|$ be the {\em underlying topological space}; the space
obtained by taking the disjoint union of the model tetrahedra and
taking the quotient by the face pairings.

At this point we should fix an encoding scheme which translates
triangulations into binary strings.  However we will not bother to do
more than remark that there are naive schemes which require about $n
\log(n)$ bits to specify a triangulation with $n$ tetrahedra.  (This
blow-up in length is due to the necessity of giving the tetrahedra
unique names.)  Thus it is only a slight abuse of language to say that
a triangulation $T$ has size $n$ when in fact its representation as a
binary string is somewhat longer.


Recall that the {\em three-sphere} is the three-manifold:

$$S^3 = \{ x \in \RR^4 \st \|x\| = 1 \}.$$

The connect sum $M \connect N$ of two connected three manifolds $M$
and $N$ is obtained by removing an open three-ball from the interior
of each of $M$ and $N$ and gluing the resulting two-sphere boundary
components.  The connect sum naturally extends to a collection of
connected three-manifolds; if $M$ is the disjoint union of connected
three-manifolds then $\connect M$ denotes their connect sum.

Note that by Alexander's Theorem (Theorem~1.1 of~\cite{Hatcher01}) $M
\connect S^3$ is homeomorphic to $M$, for any three-manifold $M$.  So
we adopt the convention that the empty connect sum yields the
three-sphere.

We now give a slightly non-standard definition of {\em compression
body}.  Take $S$ a closed orientable surface.  Let $C_0 = S \cross
[0,1]$.  Choose a disjoint collection of simple closed curves in some
component of $S \cross \{0\}$ and attach two-handles in the usual
fashion along these curves.  Cap off some (but not necessarily all) of
any resulting two-sphere boundary components with three-handles.  The
final result, $C$, is a {\em compression body}.  Set $\bdy_+ C = S
\cross \{1\}$ and set $\bdy_- C = \bdy C \setminus \bdy_+ C$.  Our
definition differs from others (\eg,~\cite{CassonGordon87}) in that
two-sphere components in $\bdy_- C$ are allowed.  The reasons for this
are explained in Remark~\ref{Rem:CompressionBodyDef}.

\section{Preliminaries}
\label{Sec:Prelim}

Here we give a few algorithms which take triangulations and check
topological properties.  See~\cite{HassEtAl99} for a more in-depth
discussion of these.


\begin{theorem}
\label{Thm:ThreeManifoldAlg}
There is a polynomial time algorithm which, given a triangulation $T$,
decides whether or not $|T|$ is a three-manifold.
\end{theorem}

\begin{proof}
Apply Poincare's Theorem: Suppose the frontier of a regular
neighborhood for each vertex of the triangulation is either a two-disk
or a two-sphere.  Then $|T|$ is a three-manifold.  The converse also
holds.

There are only linearly many (in $\size(T)$) vertices.  Checking that
each frontier is a two-sphere or a disk takes time at most polynomial
(again, in $\size(T)$) as each is a union of at most linearly many
normal triangles (see Section~\ref{Sec:NormalSurfaces}).  Thus
checking the hypothesis of Poincare's Theorem takes time at most
polynomial in $\size(T)$.
\end{proof}

Recall that a three-manifold $M$ is a {\em homology three-sphere} if
it has the same homology groups as $S^3$.  

\begin{theorem}
\label{Thm:HomologySphereAlg}
There is a polynomial time algorithm which, given a triangulation $T$
of a three-manifold, decides whether or not $|T|$ is a homology
three-sphere.
\end{theorem}

\begin{proof}
First apply Theorem~\ref{Thm:ThreeManifoldAlg} to check that $|T|$ is
in fact a three-manifold.  The homology groups $H_*(|T|, \ZZ)$ may be
read off from the Smith Normal Form of the chain boundary maps
(see~\cite{ChangDonald91}, Section 2, for an accessible overview of
algorithmic computation of homology).
Smith Normal Form of an integer matrix may be computed in polynomial
time (see~\cite{Iliopoulos89}).
\end{proof}

We also record, for future use, a few consequences of the homology
three-sphere assumption:

\begin{lemma}
\label{Lem:ConsequencesOfHomologySphere}
If $M^3$ is a homology three-sphere then $M$
is connected, closed, and orientable.  Also every closed, embedded
surface in $M$ is orientable and separating.  Finally, every connect
summand of $M$ is also a homology three-sphere.
\end{lemma}

In particular no lens space (other than $S^3$) appears as a summand of
a homology three-sphere.

We end this section with the simple:

\begin{theorem}
\label{Thm:IdenticalAlg}
There is a polynomial time algorithm which, given triangulations $T$ and
$T'$, decides whether or not $T$ is identical to $T'$.
\end{theorem}

\begin{proof}
Recall that $T$ and $T'$ are labelled: all of the tetrahedra come
equipped with names.  To check for isomorphism simply check that every
name appearing in $T$ also appears in $T'$ and that all of the face
pairings in $T$ and $T'$ agree.
\end{proof}

\begin{remark}
\label{Rem:IsomorphismOfThreeManifolds}
Note that the labelling is not needed to determine isomorphism of
triangulations.  This is because an isomorphism is completely
determined by the image of a single tetrahedron.
\end{remark}

\section{Normal and almost normal surfaces}
\label{Sec:NormalSurfaces}

In order to study triangulations we first discuss Haken's theory of
{\em normal surfaces}.  See~\cite{HassEtAl99} for a more through
treatment.

On a face $f$ of the model tetrahedron $\tau$ there are three kinds of
properly embedded arc with end points in distinct edges of $f$.  These
are called {\em normal arcs}.  A simple close curve $\alpha \subset
\bdy \tau$ is a {\em normal curve} if $\alpha$ is transverse to the
one-skeleton of $\tau$ and $\alpha$ is a union of normal arcs.  The
{\em length} of a normal curve $\alpha$ is the number of normal arcs
it contains.  A normal curve $\alpha$ is called {\em short} if it has
length three or four.

\begin{lemma}
\label{Lem:ShortCurves}
A normal curve $\alpha \subset \bdy \tau$ misses some edge of $\tau^1$
if and only if $\alpha$ meets every edge at most once if and only if
$\alpha$ is short
\end{lemma}

To see this let $\{v_{(i,j)} \st 0 \leq i < j \leq 3\}$ be the number
of intersections of $\alpha$ with each of the six edges of $\tau$.
There are twelve inequalities $v_{(0,1)} \leq v_{(1,2)} + v_{(0,2)}$,
\etc, as well as six equalities $v_{(0,1)} + v_{(1,2)} + v_{(0,2)} = 0
\mod 2$, \etc.  Easy calculation gives the desired result.

\begin{figure}
$$\begin{array}{cc}
\epsfig{file=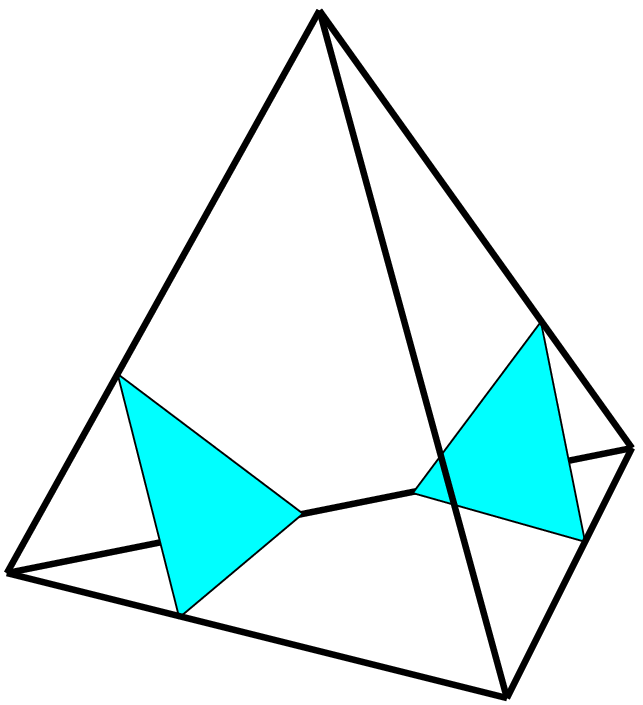, height = 3.5 cm} &
\epsfig{file=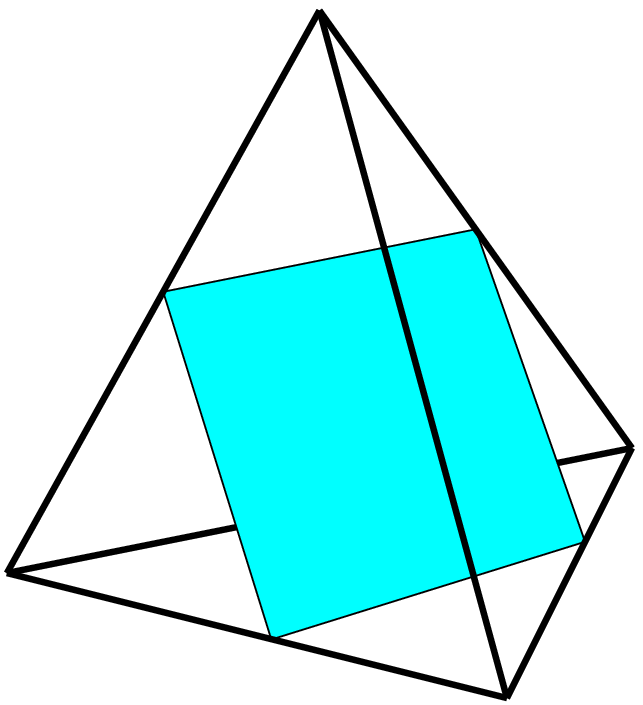, height = 3.5 cm}
\end{array}$$
\caption{Two of the four triangles and one of the three quads.}
\label{Fig:NormalDisks}
\end{figure}

In a model tetrahedron there are seven types of {\em normal disk},
corresponding to the seven distinct short normal curves in $\bdy
\tau$.  See Figure~\ref{Fig:NormalDisks}.  These are the four {\em
normal triangles} and three {\em normal quads}.  We have triangles of
type $0$, $1$, $2$, or $3$ depending on which vertex they cut off of
the model tetrahedron, $\tau$.  We have quads of type $1$, $2$, or $3$
depending on which vertex is grouped with $0$ when we cut $\tau$ along
the quad.

A surface $S$ properly embedded in $|T|$ is {\em normal} if $S \cap
\tau$ is a collection of normal disks for every tetrahedron $\tau \in
T$.

\begin{figure}
$$\begin{array}{cc}
\epsfig{file=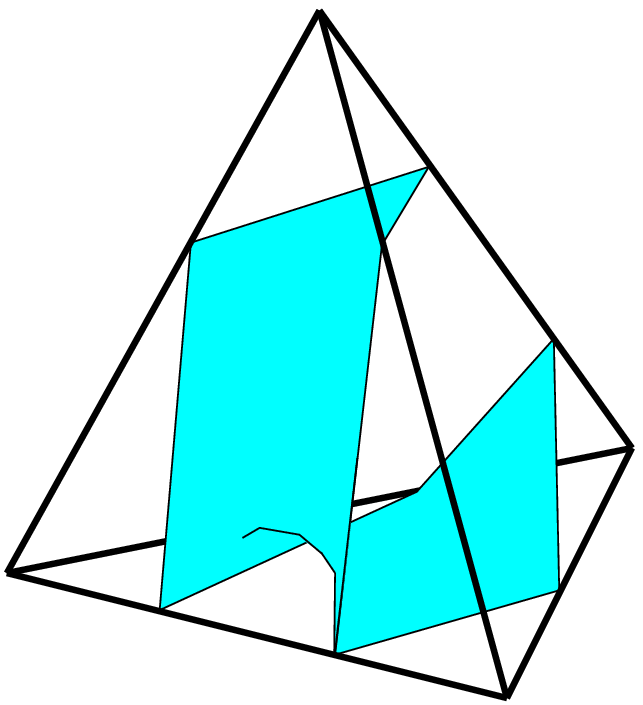, height = 3.5 cm} &
\epsfig{file=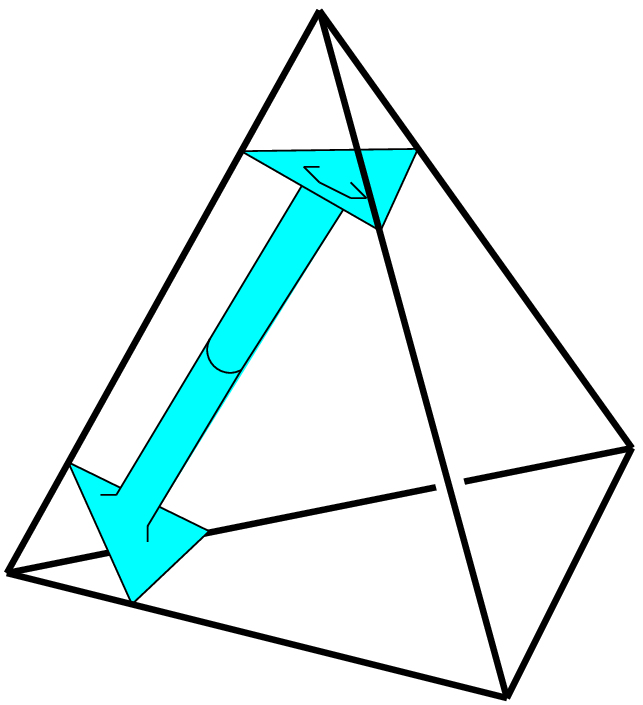, height = 3.5 cm}
\end{array}$$
\caption{One of the three octagons and one of the 25 annuli.}
\label{Fig:AlmostNormalPieces}
\end{figure}

There is also the almost normal octagon and almost normal annulus,
defined by Rubinstein~\cite{Rubinstein97}.  See
Figure~\ref{Fig:AlmostNormalPieces} for examples.  An octagon is a
disk in the model tetrahedron bounded by a normal curve of length
eight.  An annulus is obtained by taking two disjoint normal disks and
tubing them together along an arc parallel to an edge of the model
tetrahedron.  A surface $S$ properly embedded in $|T|$ is {\em almost
normal} if $S \cap \tau$ is a collection of normal disks, for every
tetrahedron $\tau \in T$, except one.  In the exceptional tetrahedron
there is a collection of normal disks and exactly one almost normal
piece.

\subsection{Weight and Euler characteristic}

For an either normal or almost normal surface $S$ take the {\em
weight} of $S$, $\weight(S) = |S \cap T^1|$, to be the number of
intersections between $S$ and the one-skeleton $T^1$.  Note that
$\size(S)$, the number of bits required to describe $S$, is about
$\size(T) \log(\weight(S))$.  To see this record a normal surface $S$
as a {\em surface vector} $v(S) \in \ZZ^{7\cdot\size(T)}$ where the
first $4\cdot\size(T)$ coordinates describe the number of normal
triangles of each type while the last $3\cdot\size(T)$ coordinates
describe the number of normal quads of each type.  At least two-thirds
of these last $3\cdot\size(T)$ coordinates are zero, as an embedded
surface has only one kind of normal quad in each tetrahedron.

For an almost normal surface $S$ we again record the vector $v(S)$ of
numbers of normal disks, as well as the type of the almost normal
piece and the name of the tetrahedron containing it.

If two normal (or almost normal) surfaces $S$ and $S'$ have the same
vector then $S$ is {\em normally isotopic} to $S'$.  This is the
natural equivalence relation on these surfaces.  As such we refer to
normal or almost normal surfaces and their vectors interchangeably
where this does not cause confusion.

We now have a few results concerning normal and almost normal
surfaces.  We assume throughout that the triangulation $T$ has
underlying space a three-manifold.

\begin{lemma}
\label{Lem:EulerAlg}
There is an polynomial time algorithm which, given a triangulation $T$
and a normal or almost surface vector $v(S)$, computes the
weight of $S$ and the Euler characteristic of $S$.
\end{lemma}

\begin{proof}
To find the weight of $S$ on a single edge $e$ of $T^1$ count
the number of normal disks meeting $e$ (with multiplicity depending on
how many times the containing tetrahedron meets $e$) and divide by the
valency of $e$ in $T^2$, the two-skeleton.

For the Euler characteristic simply use the formula $\euler(S) = F - E
+ V$ and the cell structure on $S$ coming from its being a normal
surface.  (If $S$ contains an almost normal annulus then we must add
a single edge running between the two boundary components of the
annulus.) Counting the number of faces and edges is straight-forward.
The number of vertices equals the weight.

See~\cite{AgolEtAl02}, the end of Section~5, for a more detailed
discussion.
\end{proof}

\begin{theorem}[Agol, Hass Thurston~\cite{AgolEtAl02}]
\label{Thm:ConnectedAlg}
There is an polynomial time algorithm which, given a triangulation $T$
and a normal or almost normal surface vector $v(S)$, produces the
surface vectors for the connected components of $S$.
\end{theorem}

A caveat is required here -- if several normally isotopic copies of $F$
appear in $S$ then the algorithm of Theorem~\ref{Thm:ConnectedAlg}
produces $v(F)$ only once and also reports the number of copies.  This
is required if the algorithm is to run in polynomial time on input of
the form $n\cdot F + G$.

\begin{proof}[Proof of Theorem~\ref{Thm:ConnectedAlg}]
This is one application of the ``extended counting algorithm'' given
in~\cite{AgolEtAl02}.  See the proof of Corollary~17 of that paper.
\end{proof}

\subsection{Vertex linking}

Fix a triangulation $T$ of some three-manifold.  Suppose $x \in |T|$
is a vertex of the triangulation.  Let $S$ be the frontier of small
regular neighborhood of $x$.  Then $S$ is a connected normal surface
which contains no normal quads.  Such a surface is a {\em vertex
linking} or simply a {\em vertex link}.  A normal sphere which
contains a normal quad is called {\em non-trivial}.  If the
triangulation contains no non-trivial normal two-spheres then $T$ is
{\em zero-efficient}.

\subsection{Haken sum}

Suppose $S, F, G$ are three non-empty normal surfaces with $v(S) =
v(F) + v(G)$.  Then $S$ is the {\em Haken sum} of $F$ and $G$ or,
equivalently, $S$ {\em decomposes} as a Haken sum.  Likewise, suppose
$S$ and $F$ are almost normal with identical almost normal piece, $G$
is normal, and again $v(S) = v(F) + v(G)$.  Then, again, $S$ is a
Haken sum.  In either case we write $S = F + G$.  If $S$, normal or
almost normal, does not decompose as a Haken sum then $S$ is {\em
fundamental}.  As an easy exercise:

\begin{lemma}
\label{Lem:VertexSummand}
If $S = F + G$, where $G$ is a vertex link, then $S$ is not
connected. \qed
\end{lemma}

Also,

\begin{lemma}
\label{Lem:FundamentalBound}
If $S \subset |T|$ is a fundamental normal or almost normal surface then
the largest entry of $v(S)$ is at most $\exp(\size(T))$. 
\end{lemma}

That is, there is a constant $c$ (not depending on $T$ or $S$) such
that the largest entry is less than $2^{c \cdot \size(T)}$.  This
lemma is proved for normal surfaces by Lemma~6.1 of~\cite{HassEtAl99}
and their proof is essentially unchanged for almost normal
surfaces. \qed

\begin{lemma}
\label{Lem:FundamentalTwoSphere}
Suppose $T$ is a triangulation of a homology three-sphere.  Suppose
$T$ contains a non-trivial normal two-sphere.  Then $T$ contains a
non-trivial normal two-sphere which is fundamental.
\end{lemma}

\begin{proof}
This is discussed as Proposition~4.7 of~\cite{JacoRubinstein02a}.  The
essential points are that Euler characteristic is additive under Haken
sum, that $T$ does not contain any normal $\RRPP^2$ or $\DD^2$ (by
Lemma~\ref{Lem:ConsequencesOfHomologySphere}), and that no
summand is vertex-linking (by Lemma~\ref{Lem:VertexSummand}).
\end{proof}

\begin{lemma}
\label{Lem:FundamentalAlmostTwoSphere}
Suppose $T$ is zero-efficient triangulation of a homology
three-sphere.  Suppose $T$ contains an almost normal two-sphere.  Then
$T$ contains a fundamental almost normal two-sphere.
\end{lemma}

\begin{proof}
This is identical to the proof of
Lemma~\ref{Lem:FundamentalTwoSphere}, except that $S$ cannot have any
normal two-sphere summand as $T$ is zero-efficient.
\end{proof}

Of a much different level of difficulty is Thompson's Theorem:

\begin{theorem}[Thompson~\cite{Thompson94}]
\label{Thm:ThompsonSphere}
Suppose $|T| \homeo S^3$.  Suppose also that $T$ is zero-efficient.
Then $T$ contains an almost normal two-sphere. \qed
\end{theorem}

We cannot do better than refer the reader to Thompson's original paper
and remark that the proof uses Gabai's notion of {\em thin position}
for knots~\cite{Gabai87}.

We end this section with:

\begin{lemma}
\label{Lem:FundamentalSphereAlg}
There is an exponential time algorithm which, given a triangulation
$T$ of a three-manifold, produces either the surface vector of a
fundamental non-trivial normal two-sphere or, if none exists, produces
the surface vector of a fundamental almost normal two-sphere or, if
neither exists, reports ``$|T|$ is not homeomorphic to the
three-sphere''.
\end{lemma}

We only sketch a proof -- the interested reader should
consult~\cite{HassEtAl99}, \cite{JacoRubinstein02a} (page 53),
or~\cite{Barchechat03} (page 83).  If $T$ is not zero-efficient there
is a fundamental normal two-sphere.  This can be found by enumerating
all fundamental surfaces (a finite list, by work of Haken) and
checking each surface on the list.  If $T$ is zero-efficient then no
non-trivial normal two-sphere appears.  However we again have that
some fundamental almost normal two-sphere exists, if $T$ contains any
almost normal two-sphere.  Finally, if no non-trivial normal sphere
nor any almost normal sphere appears in amongst the fundamentals then,
combining Lemma~\ref{Lem:FundamentalAlmostTwoSphere} and
Theorem~\ref{Thm:ThompsonSphere}, conclude that $|T|$ is not the
three-sphere.

As presented the running time of the algorithm is unclear; it depends
on the number of fundamental surfaces.  However, using vertex
fundamental surfaces (see Remark~\ref{Rem:VertexFundamental}) and
linear programming techniques Casson reduces the search to take time
at most a polynomial times $3^{\size(T)}$. \qed

\section{Blocked submanifolds}
\label{Sec:Blocked}

Normal (and almost normal) surfaces cut a triangulated manifold into
pieces.  These submanifolds have natural polyhedral structures which
we now investigate.

\begin{figure}
$$\begin{array}{c}
\epsfig{file=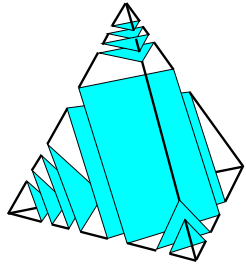, height = 3.5 cm}
\end{array}$$
\caption{The tetrahedron $\tau$ cut along $S$.  Note that there are
  two blocks of the form ``normal disk cross interval''.}
\label{Fig:Blocks}
\end{figure}

Let $\tau$ be a model tetrahedron, and suppose that $S \subset \tau$
is a embedded collection of normal disks and at most one almost normal
piece.  Let $B$ be the closure of any component of $\tau \setminus
S$.  We call $B$ a {\em block}.  See Figure~\ref{Fig:Blocks}.

An block containing exactly two normal disks of the same type is
called a {\em product block}.  All other blocks are called {\em core
blocks}.  Note that there are only seven kinds of product block
possible, corresponding to the seven types of normal disks.  Likewise
there is a bounded number of core blocks.  Five such are shown in
Figure~\ref{Fig:Blocks}, but many more are possible.  Most of these
meet an almost normal annulus.

The components of $\bdy B$ meeting $S$ are the {\em horizontal}
boundary components of $B$, denoted $\hbdy B$. All other faces of $B$
(the faces of $B$ which lie in the two-skeleton) are $\vbdy B$, the
{\em vertical} boundary.

Suppose now that $T$ is a triangulation of a three-manifold
and $S \subset |T|$ is a normal or almost normal surface.  For
simplicity, suppose that $S$ is a transversely oriented and
separating.  Let $N_S$ be the closure of the component of $|T|
\setminus S$ into which the transverse orientation points.

Then $N_S$ is a {\em blocked submanifold} of $|T|$.  Note that the
induced cell structure on $N_S$ (coming from $T$) naturally breaks into
two parts.  So, let $\blocky{N}_P$ be the union of all product blocks
in $N_S$ and let $\blocky{N}_C$ be the union of all core blocks in $N_S$.

\begin{remark}
\label{Rem:BoundedCore}
In any blocked submanifold there are at most a linear number (in
$\size(T)$) of core blocks.  In fact there at most six in each
tetrahedron (plus possibly two more coming from the almost normal
annulus).  This simple observation underlies Kneser-Haken finiteness
(see, for example,~\cite{Haken68}) and is generally useful in the
algorithmic setting.
\end{remark}

Note that $\blocky{N}_P$ and $\blocky{N}_C$ are not necessarily
submanifolds of $|T|$.  To produce submanifolds take $N_P$ to be a
regular closed neighborhood of $\blocky{N}_P$, where the neighborhood
is taken inside of $N_S$.  Also, take $N_C$ to be the closure of $N_S
\setminus N_P$.  Note the asymmetry between the definitions of $N_P$
and $N_C$: we have $\blocky{N}_P \subset N_P$ while $N_C \subset
\blocky{N}_C$.  As above define $\hbdy N_P = N_P \cap S$ and $\vbdy
N_P = \closure{\bdy N_P \setminus \hbdy N_P}$.  The horizontal and
vertical boundaries $\hbdy N_C$ and $\vbdy N_C$ are defined similarly.

Note also that $\blocky{N}_P$ may be represented by a {\em block
vector}; an element $v(\blocky{N}_P) \in \ZZ^{7\cdot\size(T)}$ where
the first $4\cdot\size(T)$ coordinates describe the number of triangle
product blocks of each type while the last $3\cdot\size(T)$
coordinates describe the number of quad product blocks of each type.

We now have:

\begin{theorem}
\label{Thm:AHTBlockVectorAlg}
There is an polynomial time algorithm which, given a triangulation $T$
and blocked submanifold $N_S \subset |T|$ (via the surface vector
$v(S)$), produces the block vector for each connected component of
$\blocky{N}_P$. 
\end{theorem}

We only require this algorithm for blocked submanifolds $N_S$ cut out
of $|T|$ by a transversely oriented separating normal or almost normal
two-sphere $S$.  In this case, the input to the algorithm is just $T$
and the surface vector $v(S)$.  From $v(S)$ it is possible (in time
polynomial in $\size(T)$ and $\log(\weight(S))$) to find a block
vector for $\blocky{N}_P$.  Given this the other details for a proof
of Theorem~\ref{Thm:AHTBlockVectorAlg} can be found
in~\cite{AgolEtAl02}. \qed

\begin{remark}
\label{Rem:NumberOfN_PComponents}
We also remark that $\blocky{N}_P$ has at most a linear number (in
$\size(T)$) of connected components.  This is because $\vbdy N_P =
\vbdy N_C$ and the latter has at most linearly many components.  (See
Remark~\ref{Rem:BoundedCore}.)  Thus, unlike the algorithm of
Theorem~\ref{Thm:ConnectedAlg} we do not need to concern ourselves
with components of $\blocky{N}_P$ appearing with high multiplicity.
\end{remark}

\section{Normalizing slowly}
\label{Sec:OldNormalization}

In this section we discuss a restricted version of Haken's {\em
normalization} procedure for producing normal surfaces.  This material
appeared first in an unpublished preprint of mine and later in my
thesis~\cite{Schleimer01}.  I thank Danny Calegari for reading an
early version of this work.  I also thank Bus Jaco for several
enlightening conversations regarding Corollary~\ref{Cor:Barrier}.

Several authors have independently produced versions of these ideas;
for example see~\cite{JacoRubinstein02a}, \cite{Barchechat03}, or
\cite{King01}.

Let $S \subset |T|$ be a transversely oriented, almost normal surface.
Here $T$ is triangulation of a closed, orientable, connected
three-manifold.

\begin{define}
A compression body $C_S \subset |T|$ is {\em canonical} for $S$ if
$\bdy_+ C_S = S$, $\bdy_- C_S$ is normal, the transverse orientation
points into $C_S$, and any normal surface $S'$ disjoint from $S$ may
be normally isotoped to one disjoint from $C_S$.
\end{define}

As a bit of notation take $\bdy_- C_S = \normal{S}$ and call
this the {\em normalization} of $S$.

\vspace{2mm}
\noindent
{\bf Theorem~\ref{Thm:NormalizationAlg}.}
{\em 
Given a transversely oriented almost normal surface $S$ there exists a
canonical compression body $C_S$ and $C_S$ is unique (up to normal
isotopy).  Furthermore there is a algorithm which, given the
triangulation $T$ and the surface vector $v(S)$, computes the surface
vector of $\bdy_- C_S = \normal{S}$.
}
\vspace{1mm}

The proof of this theorem is lengthy and is accordingly spread from
Section~\ref{Sec:Tight} to Section~\ref{Sec:NormalizationAlg}.  We
here give the necessary definitions.  In Section~\ref{Sec:Tight} we
discuss the tightening procedure.  In Section~\ref{Sec:Embedded} we
show that the tightening procedure gives an embedded isotopy.  We
discuss the capping off procedure in Section~\ref{Sec:CappingOff}.
Finally in Section~\ref{Sec:NormalizationAlg} we prove
Theorem~\ref{Thm:NormalizationAlg}.

\subsection{Non-normal surfaces}

Let $S$ be a surface properly embedded in a triangulated
three-manifold $|T|$ and suppose that $S$ is transverse to the skeleta
of $T$.  Denote the $i$-skeleton of $T$ by $T^i$.  Generalize the
notion of weight so that $\weight(S) = |S \cap T^1|$.

\begin{figure}
\begin{center}
\epsfig{file=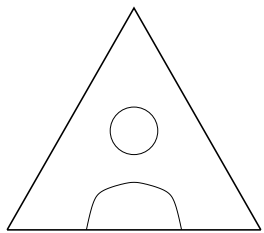, height=3cm}
\end{center}
\caption{A simple curve and a bent arc.}
\label{Fig:NonNormalCurves}
\end{figure}

We characterize some of the ways $S$ can fail to be normal.  A {\em
simple curve} of $S$ is a simple closed curve of intersection between
$S$ and the interior of some triangular face $f \in T^2$.  Also, a
{\em bent arc} of $S$ is a properly embedded arc of intersection
between $S$ and the interior of some triangular face $f \in T^2$ with
both endpoints of the arc contained in a single edge of $f$. Both of
these are drawn in Figure~\ref{Fig:NonNormalCurves}.

\subsection{Surgery and tightening disks}

\begin{define}
An embedded disk $D \subset |T|$ is a {\em surgery disk} for $S$ if $D
\cap S = \bdy D$, $D \subset T^2$ or $D \cap T^2 = \emptyset$, and
$D \cap T^1 = \emptyset$.
\end{define}

There is a {\em surgery} of $S$ along $D$: Remove a small neighborhood
of $\bdy D$ from $S$ and cap off the boundaries thus created with
disjoint, parallel copies of $D$.  Note that we do not require $\bdy
D$ to be essential in $S$.  A simple curve of $S \cap T^2$ is {\em
innermost} if it is the boundary of a surgery disk embedded in a
triangle of $T^2$.

\begin{define}
An embedded disk $D \subset |T|$ is a {\em tightening disk} for $S$ if
$\bdy D = \alpha \cup \beta$ where $\alpha \cap \beta = \bdy \alpha =
\bdy \beta$ and $D \cap S = \alpha$.  Also, either $D \subset T^2$ or
$D \cap T^2 = \beta$.  But in either case $D \cap T^1 = \beta$, and $D
\cap T^0 = \emptyset$.
\end{define}

There is a {\em tightening isotopy} of $S$ across $D$: Push $\alpha$
along the disk $D$, via ambient isotopy of $S$ supported in a small
neighborhood of $D$, until $\alpha$ moves past $\beta$. This procedure
reduces $\weight(S)$ by exactly two. A bent arc of $S$ is {\em
outermost} if it lies on the boundary of a tightening disk embedded in
a triangle of $T^2$.

Suppose $S$ contains an almost normal octagon.  There are two
tightening disks on opposite sides of the octagon both giving
tightening isotopies of $S$ to a possibly non-normal surface of lesser
weight.  (See Figure~\ref{Fig:AlmostNormalPieces}.)  We call these the
{\em exceptional tightening disks}.  If $S$ contains an almost normal
annulus then the tube is parallel to at least one edge of the
containing tetrahedron.  (Again see
Figure~\ref{Fig:AlmostNormalPieces}.)  For every such edge there is an
exceptional tightening disk.  Also, the disk which surgers the almost
normal annulus will be called the {\em exceptional surgery disk}.

\subsection{Normal isotopy}

We sharpen our notion of equivalence:

\begin{define}
An isotopy $H: |T| \times I \rightarrow |T|$ is a {\em normal isotopy}
if, for all $s \in I$ and for every simplex $\sigma$ in $T$,
$H_s(\sigma) = \sigma$.
\end{define}

Two subsets of $|T|$ are {\em normally isotopic} if there is a normal
isotopy taking one to the other.

\section{Tightening}
\label{Sec:Tight}

This section discusses the {\em tightening} procedure which will yield
an embedded isotopy.  This is proved in Lemma~\ref{Lem:Embedded}
below.

Suppose that $S \subset |T|$ is a transversely orientable separating
almost normal surface.  Here $T$ is a triangulation of a
three-manifold.  We wish to isotope $S$ off of itself while
continuously reducing the weight of $S$.

Suppose that $D$ is an exceptional tightening disk for $S$.  Choose
the transverse orientation for $S$ which points into the component of
$|T| \setminus S$ which meets $D$.  The {\em $F$-tightening
procedure} constructs a map $\calF \from S \cross [0,n] \to |T|$ as
follows:

\begin{enumerate}
\item
Let $F_0 = S$.  Take $\calF_0 \from S \cross \{0\} \to |T|$ to be
projection to the first factor.  Let $D_0 = D$.
\item
Now do a small normal isotopy of $F_0$ in the transverse direction
while tightening $F_0$ along $D_0$.  This extends $\calF_0$ to a map
$\calF_1 \from S \cross [0, 1] \to |T|$, with $F_t = \calF_1(S \cross
\{t\})$.  Note that the surface $F_1$ inherits a transverse
orientation from $F_0$.  Arrange matters so that $F_{\frac{1}{2}}$ is
the only level which is not transverse to $T^2$.  Furthermore
$F_{\frac{1}{2}}$ only has a single tangency with $T^1$ and this
tangency occurs in the middle of $\bdy D_0 \cap T^1$.
\item
At stage $k \geq 1$ there are two possibilities.  Suppose first that
$F_k$ has an outermost bent arc $\alpha$ with the transverse
orientation of $F_k$ pointing into the tightening disk $D_k$, which is
cut out of $T^2$ by $\alpha$.  Then extend $\calF_k$ to $\calF_{k+1}
\from S \cross [0, k+1] \to |T|$ by doing a small normal isotopy of
$F_k$ in the transverse direction while tightening $F_k$ across $D_k$,
the $k^\th$ tightening disk.  So $\calF_k = \calF_{k+1}|S \cross
[0,k]$ and $F_{t} = \calF_{k+1}(S \cross \{t\})$.  Note that the
surface $F_{k+1}$ inherits a transverse orientation from $F_k$.
Arrange matters so that $F_{k + \frac{1}{2}}$ is the $k + 1^\th$ level
which is not transverse to $T^2$.  Furthermore $F_{k + \frac{1}{2}}$
only has a single tangency with $T^1$ and this tangency occurs in the
middle of $\bdy D_k \cap T^1$.

If there is no outermost bent arc $\alpha \subset F_k$ then set $n =
k$ and the procedure halts.
\end{enumerate}

\begin{remark}
\label{Rem:CalFNotUnique}
As $\weight(F_{k + 1}) = \weight(F_{k}) - 2$ this process terminates.
Note also that $\calF_n$ is far from unique -- at any stage in the
procedure there may be many tightening disks to choose from.
\end{remark}

We will show in Lemma~\ref{Lem:Embedded} that the map $\calF_n \from S
\cross [0,n] \to M$ is an embedding.  Note that, by construction, $S =
F_0 = \calF_n(S \cross \{0\})$ and in general $F_t = \calF_n(S \cross
\{t\})$.  To simplify notation set $\calF = \calF_n$.

\section{Tracking the isotopy}
\label{Sec:Embedded}

In this section we analyze how $\image(\calF)$ intersects the skeleta
of the triangulation.  Let $S \subset |T|$, $\calF$, $\calF_k$, and
$F_t$ be as defined in Section~\ref{Sec:Tight}.

\begin{figure}
\psfrag{F}{$F_k$}
\psfrag{S}{$S = F_0$}
$$\begin{array}{cc}
\epsfig{file=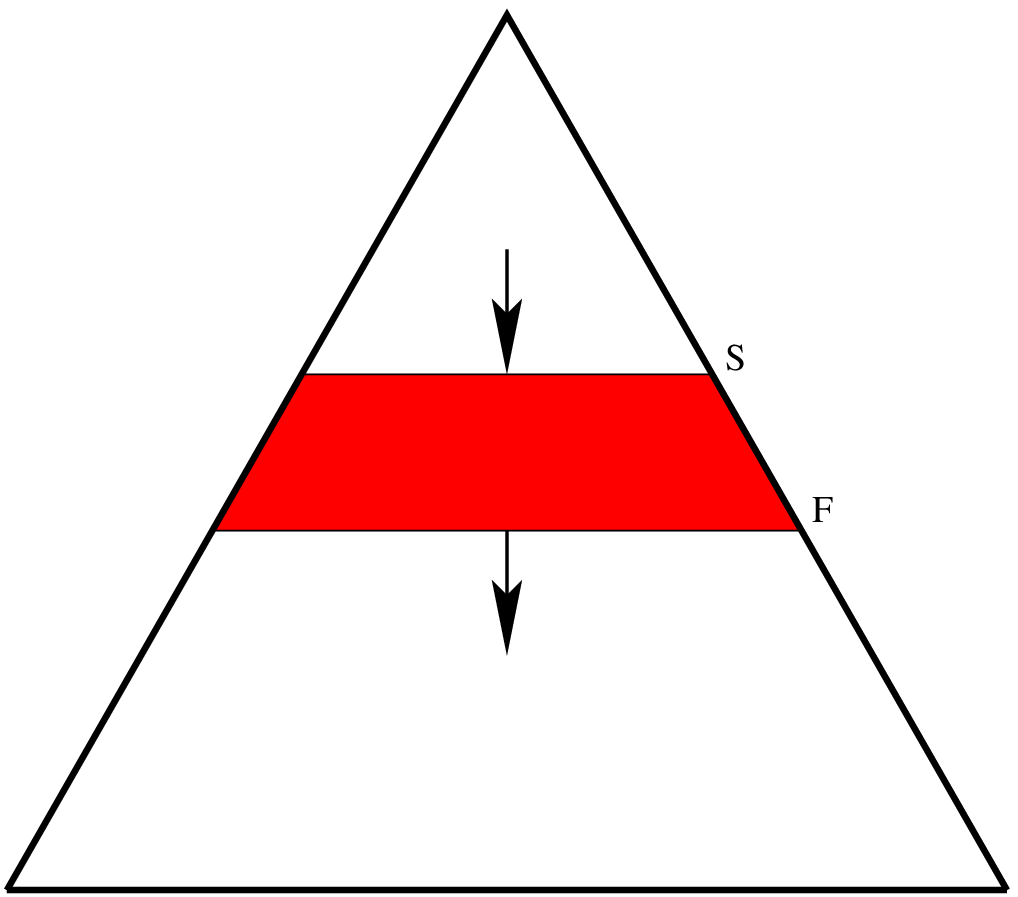, height=3cm} &
\epsfig{file=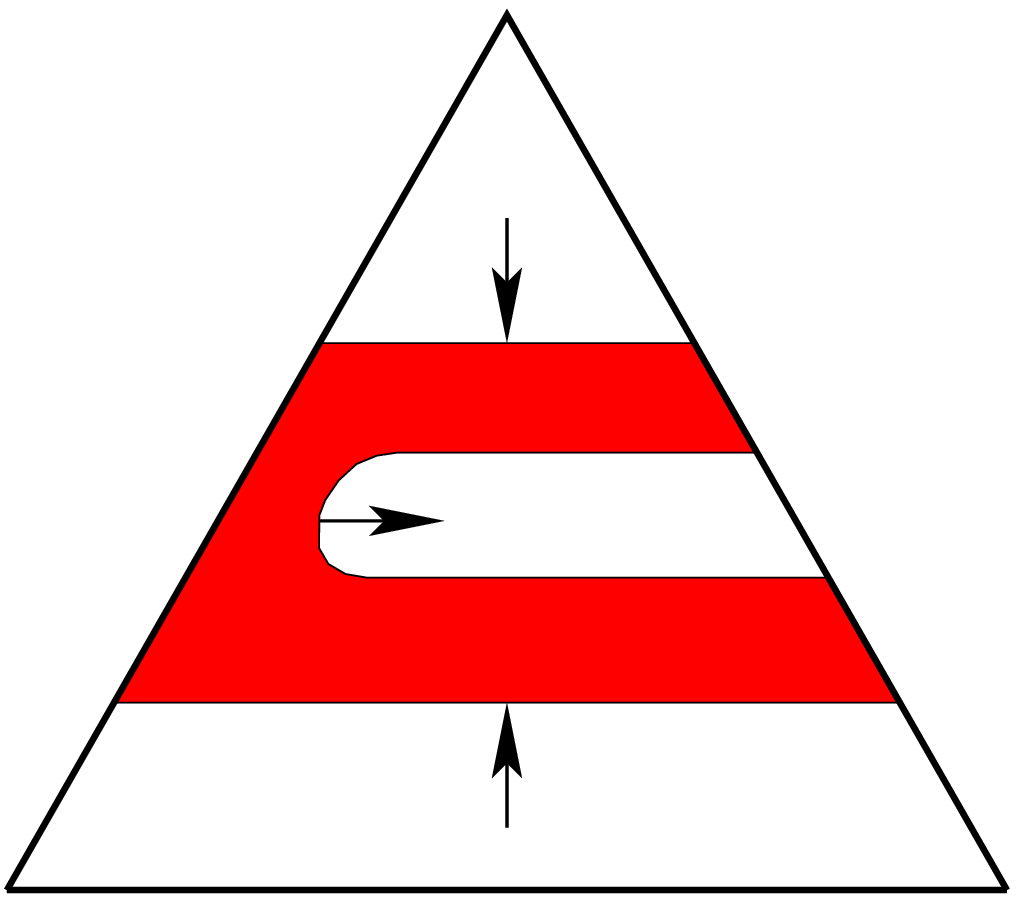, height=3cm}         \\
\mbox{Critical} & \mbox{Temporary}                      \\
\epsfig{file=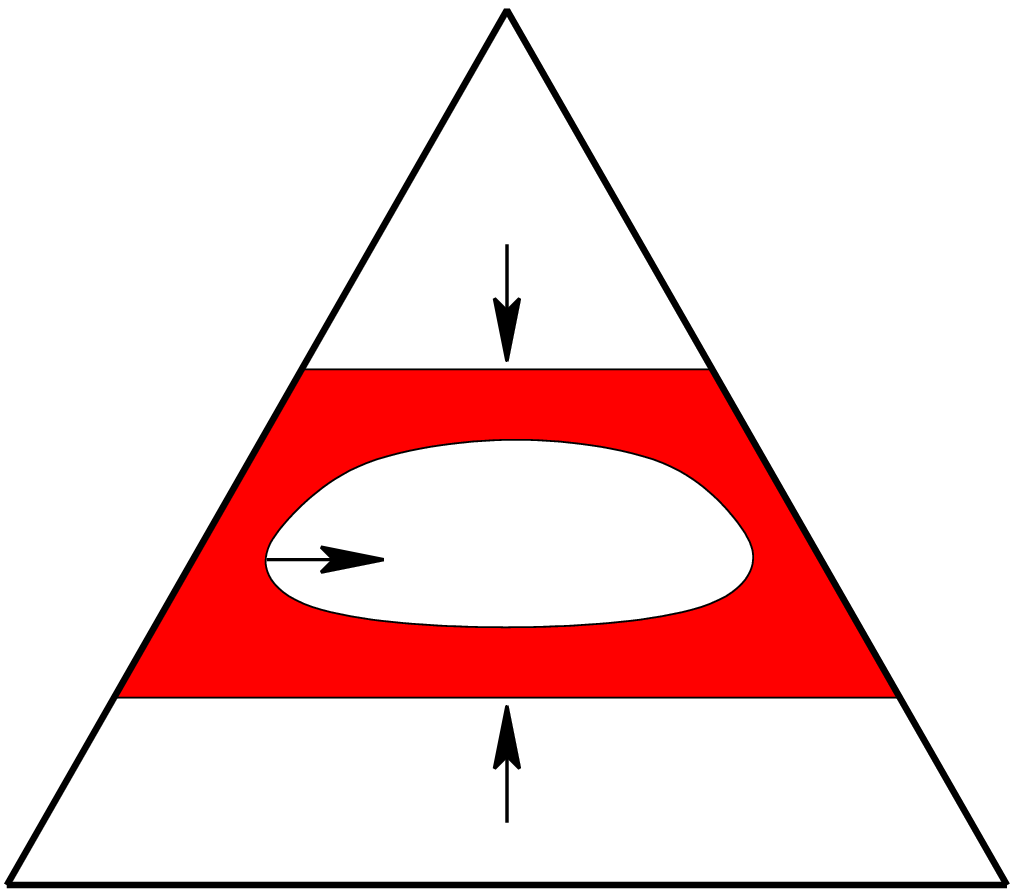, height=3cm} &
\epsfig{file=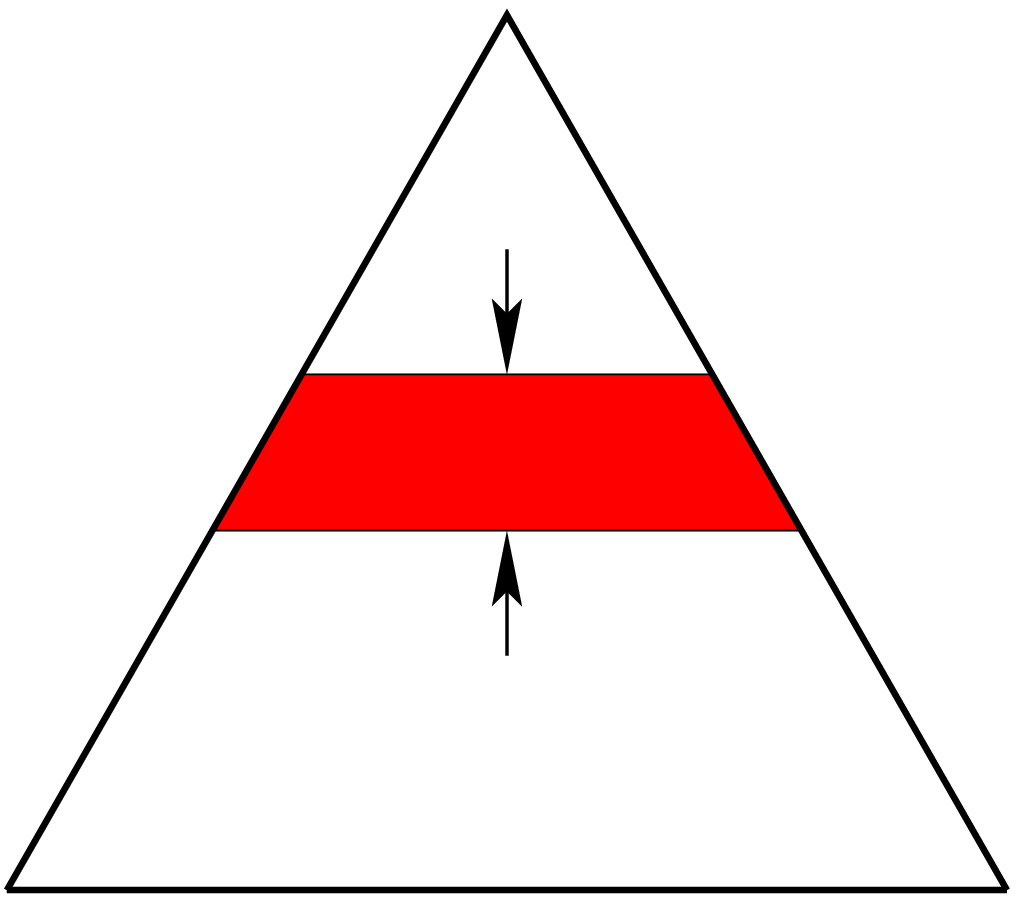, height=3cm}       \\
\mbox{Terminal with hole} & \mbox{Terminal}
\end{array}$$
\caption{The Rectangles}
\label{Fig:Rectangles}
\end{figure}

\begin{figure}
$$\begin{array}{cc}
\epsfig{file=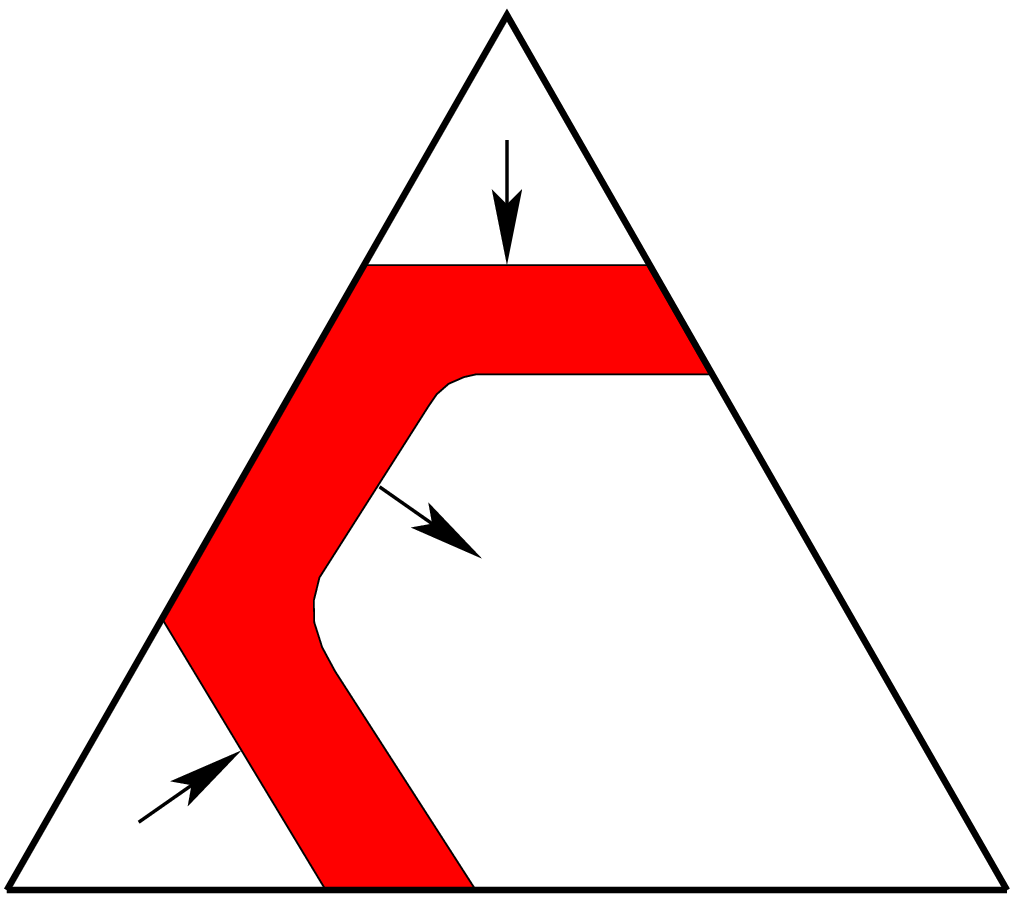, height=3cm} &
\epsfig{file=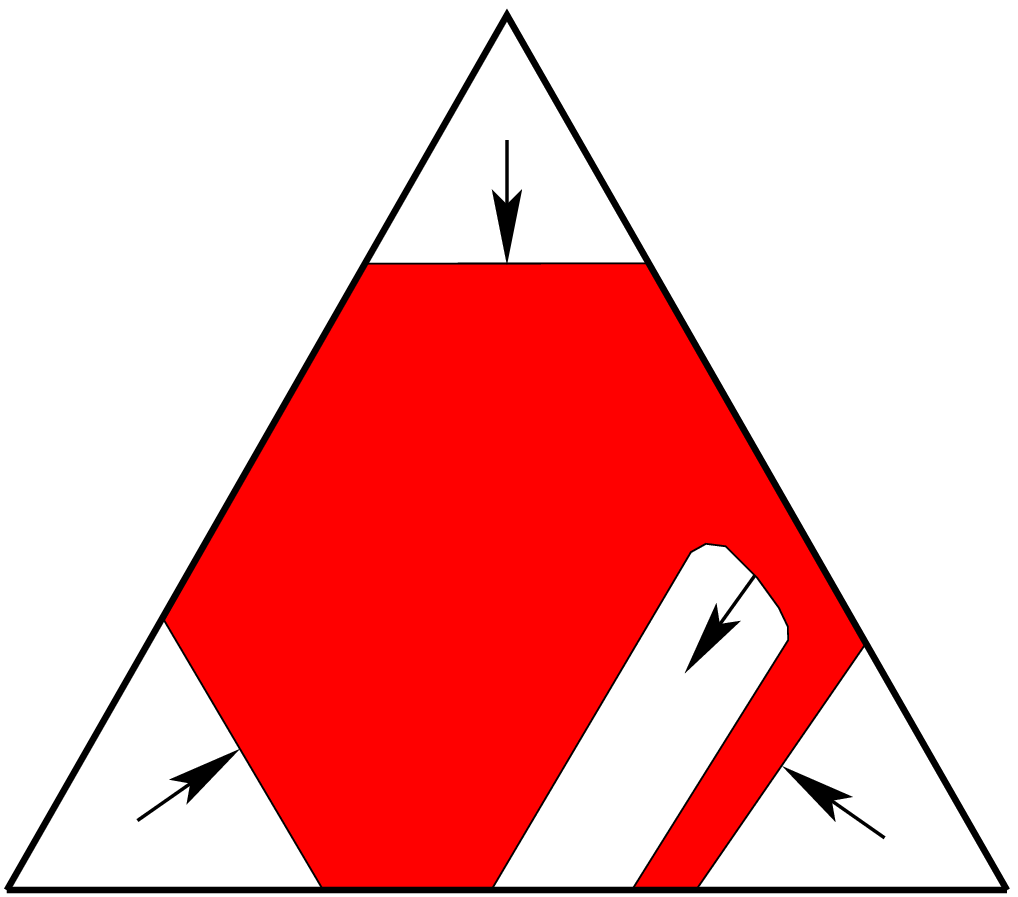, height=3cm}           \\
\mbox{Critical} & \mbox{Temporary}                      \\
\epsfig{file=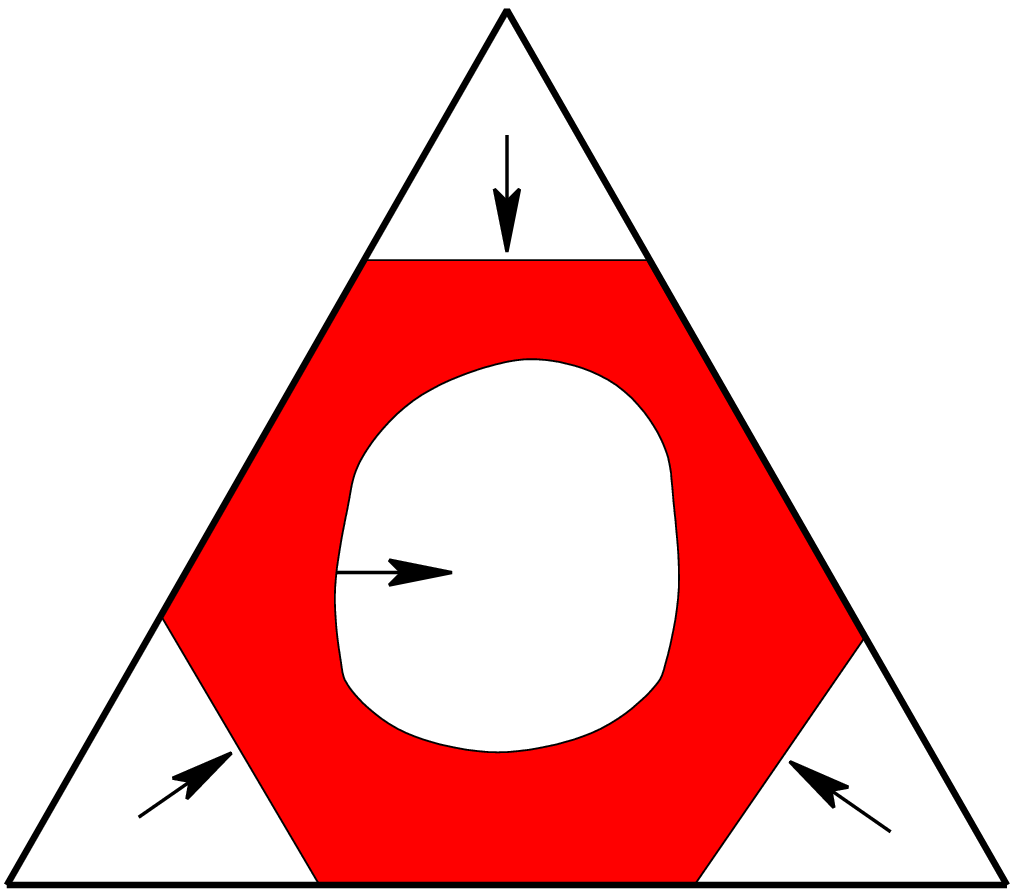, height=3cm} &
\epsfig{file=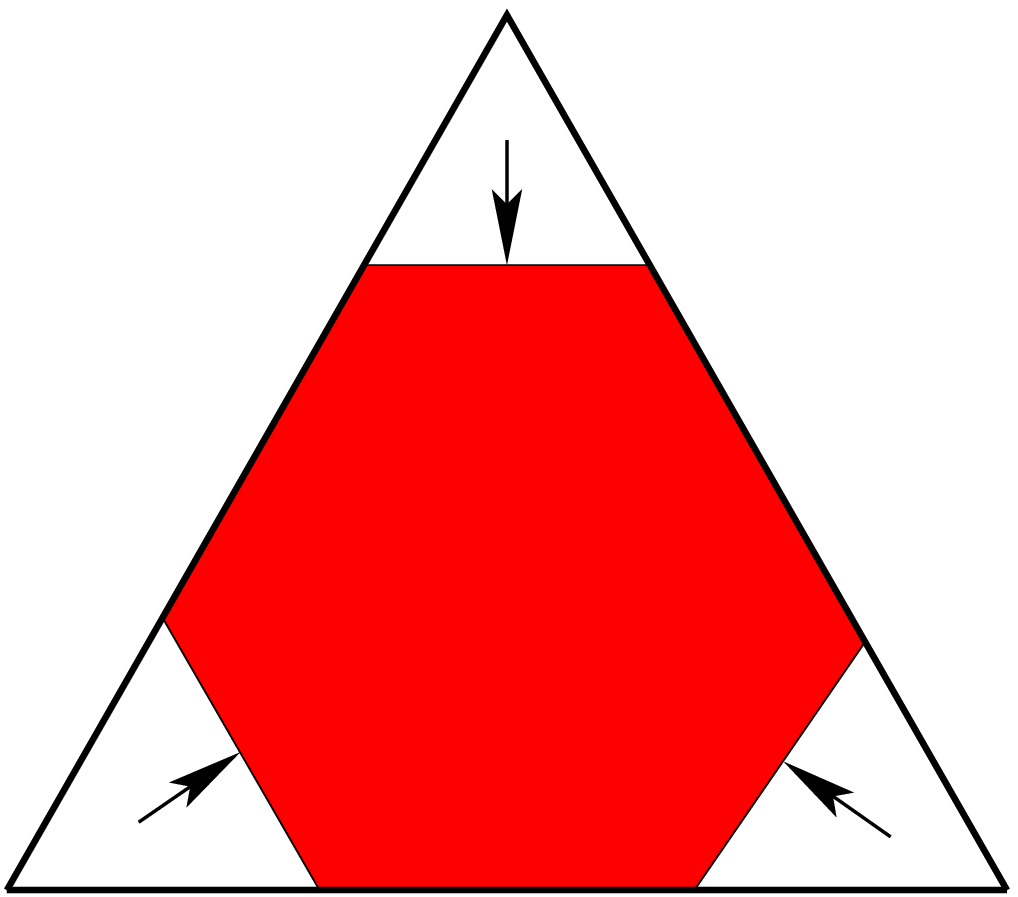, height=3cm}         \\
\mbox{Terminal with hole} & \mbox{Terminal}
\end{array}$$
\caption{The Hexagons}
\label{Fig:Hexagons}
\end{figure}

Figures~\ref{Fig:Rectangles} and~\ref{Fig:Hexagons} display a few of
the possible components of intersection $f \cap \image(\calF_k)$
assuming that $\calF_k$ is an embedding.  Here $f$ is a face of $T^2$.
Lemma~\ref{Lem:Embedded} below shows that this collection is complete
up to symmetry.  Note that the normal arcs, bent arcs, and simple
curves bounding the components receive a transverse orientation from
$S$ or $F_k$.  In these figures all arcs of $S$ are pointed towards
while arcs of $F_k$ are pointed away from, agreeing with the
transverse orientation.  The components of intersection containing a
normal arc of $F_k$ are called {\em critical}.  Those with a bent arc
of $F_k$ are called {\em temporary}.  Any component containing a
simple curve of $F_k$ is called {\em terminal with hole}.  Finally,
components of $f \cap \image(\calF_k)$ which are completely disjoint
from $F_k$ are simply called {\em terminal}.  Again, refer to
Figure~\ref{Fig:Rectangles} and~\ref{Fig:Hexagons}.

The tightening procedure combines the critical components in various
ways.  However, a temporary component always results in a terminal
(possibly with hole) and these are stable.  Note also that there
is a second critical rectangle which ``points upward.''  The
non-critical components may be foliated by the levels of $\calF_k$ in
multiple ways, depending on the order of the tightening isotopies.

\begin{lemma}
\label{Lem:Embedded}
For every $k$, the map $\calF_k$ is an embedding.  Furthermore, for $k
> 0$ and for every $f \in T^2$, the connected components of $f \cap
\image(\calF_k)$ are given, up to symmetry, by
Figures~\ref{Fig:Rectangles} and~\ref{Fig:Hexagons}.
\end{lemma}

\begin{proof}
Proceed by induction: Both claims are trivial for $k = 0$.  Now to
deal with $k = 1$.  The exceptional tightening disk $D_0$ has interior
disjoint from $S = F_0$.  It follows that $\calF_1$ is an embedding.
To verify the second claim for $k = 1$ note that the image of
$\calF_1|S \cross [0, \epsilon]$ intersects all faces $f \in T^2$ only
in critical rectangles.  Up to $t = \frac{1}{2}$ the image of
$\calF_1|S \cross [0, t]$ intersected with $f$ is combinatorially
constant.  Crossing $t = \frac{1}{2}$ adds a regular neighborhood of
$D_0$ to the image.  This only intersects $f$ in a regular
neighborhood of $\bdy D_0 \cap T^1$.  So the pieces of $f \cap
\image(\calF_1)$ are unions of critical rectangles connected by small
neighborhoods of sub-arcs of $T^1$.  Also these sub-arcs only meet the
$F_t$ side of the critical rectangles.  As each critical rectangle
meets two edges of the face $f$ it follows that at most three critical
rectangles are joined together to form a component of $f \cap
\image(\calF_1)$.  We list all possible cases -- consulting
Figures~\ref{Fig:Rectangles} and~\ref{Fig:Hexagons} will be helpful:
\begin{enumerate}
\item
Two critical rectangles in $f$ may be combined to produce a temporary
rectangle, a terminal rectangle with a hole, or a critical hexagon.
\item
Three critical rectangles in $f$ may be combined to produce a
temporary hexagon or a terminal hexagon with a hole.
\end{enumerate}

Now to deal with the general case: Suppose that both claims hold at
stage $k$.  Suppose that $\alpha \subset F_k$ is the bent arc on the
boundary of $D_k \subset f \in T^2$, the next tightening disk in the
sequence.  Suppose that $\interior(D_k)$ meets $\image(\calF_k)$.  By
the second induction hypothesis there is a component, $C$, of $f \cap
\image(\calF_k)$ which meets $\interior(D_k)$ and appears among those
listed in Figures~\ref{Fig:Rectangles} and~\ref{Fig:Hexagons}. Observe
that each component of $f \cap \image(\calF_k)$, and hence $C$, meets
at least two edges of $f$.  The bent arc $\alpha$ meets only one edge
of $f$.  It follows that the interior of $C$ must meet $\alpha$.  Thus
$\calF_k$ was not an embedding, a contradiction.

It follows that $D_k \cap \image(\calF_k) = \alpha$.  Since the
$k + 1^\th$ stage of the isotopy is supported in a small
neighborhood of $F_k \cup D_k$ it follows that $\calF_{k+1}$ is an
embedding.

Now, the transverse orientation on $F_k$ gives rise to a transverse
orientation on $F_{k+1}$.  To verify the second claim again list the
possible cases:
\begin{enumerate}
\item
Two critical rectangles in $f$ may be combined to produce a temporary
rectangle, a terminal rectangle with a hole, or a critical hexagon.
\item
Three critical rectangles $f$ may be combined to produce a temporary
hexagon or a terminal hexagon with a hole.
\item
A critical rectangle and critical hexagon in $f$ may be combined to
produce a temporary hexagon or a terminal hexagon with a hole.
\item 
A temporary component can lead only to a terminal (possibly with
hole).
\end{enumerate}
This completes the induction.
\end{proof}

\begin{remark}
\label{Rem:ShapeOfFn}
By maximality of $\calF$, the surface $F_n = \calF(S \cross \{n\})$
has no outermost bent arcs with outward orientation. A bent arc with
inward orientation would violate the second induction hypothesis of
Lemma~\ref{Lem:Embedded}. So $F_n$ contains no bent arcs. $F_n$ may
contain simple curves, but the second induction hypotheses shows that
all of these are innermost with transverse orientation pointing toward
the bounded surgery disk.
\end{remark}

Given that $\calF$ is an embedding, in the sequel $\image(\calF_k)$ is
denoted by $\calF_k$.  Replacing $S$ in Lemma~\ref{Lem:Embedded} by a
disjoint union of $S$ with a collection of normal surfaces gives:

\begin{corollary}
\label{Cor:Barrier}
If $S'$ is any normal surface in $|T|$ which does not intersect $S$
then $\calF \cap S' = \emptyset$, perhaps after a normal isotopy of
$S'$ (rel $S$).
\end{corollary}

Let $\tau$ be any tetrahedron in the given triangulation $T$.

\begin{lemma}
\label{Lem:Balls}
For all $k \geq 1$, $\tau \setminus \calF_k$ is a disjoint collection
of balls.
\end{lemma}

\begin{proof}
Again we use induction.  Our induction hypothesis is as follows: $\tau
\setminus \calF_k$ is a disjoint collection of balls, unless $k = 0$,
and $\tau$ contains the almost normal annulus of $S$.  (In this
situation $\tau \setminus \calF_0$ is a disjoint collection of balls
and one solid torus $\DD^2 \cross S^1$.)

The base case is trivial.  Suppose $B$ is a component of $\tau
\setminus \calF_k$.  There are now two cases to consider.  Either $B$
is cut by an exceptional tightening disk or it is not.  Assume the
latter.  Then $B$ is a three-ball by induction and after the $k +
1^\th$ stage of the isotopy $B \cap \calF_{k+1}$ is a regular
neighborhood (in $B$) of a collection of disjoint arcs and disks in
$\bdy B$.  Hence $B \setminus \calF_{k+1}$ is still a ball.

If $B$ is adjacent to the almost normal piece of $F_0$ then let $D_0$
be the exceptional tightening disk. Set $B_{\epsilon} = B \setminus
\neigh(D_0)$. Each component of $B_{\epsilon}$ is a ball, and the
argument of the above paragraph shows that they persist in the
complement of $\calF_1$.
\end{proof}

A similar induction argument proves:
\begin{lemma}
\label{Lem:Handlebody}
For all $k \geq 1$, $\tau \cap \calF_k$ is a disjoint collection of
handlebodies.
\end{lemma}

This lemma is not used in what follows and its proof is accordingly
left to the interested reader.  Recall that $\bdy \calF_k = S \cup
F_k$.  A trivial corollary of Lemma~\ref{Lem:Balls} is:

\begin{corollary}
\label{Cor:PlanarSurfaces}
For all $k$, the connected components of $\tau \cap F_k$ are planar.
\end{corollary}

The connected components of $\tau \cap F_n$ warrant closer attention:

\begin{lemma}
\label{Lem:BoundaryFn}
Each component of $\tau \cap F_n$ has at most one normal curve
boundary component.  This normal curve must be short.
\end{lemma}

\begin{proof}
Let $\tau \in T$ be a tetrahedron. Let $P$ be a connected component of
$\tau \cap F_n$. By Lemma~\ref{Lem:Embedded} the boundary $\bdy P$ is
a collection of simple curves and normal curves in $\bdy \tau$.  Let
$\alpha$ be any normal curve in $\bdy P$. Let $\{ \alpha_j \}$ be the
normal arcs of $\alpha$.

\begin{proofclaim}
$\alpha$ has length three or four.
\end{proofclaim}

Call the collection of critical rectangles and hexagons (in $\bdy \tau
\cap \calF$) meeting $\alpha$ the {\em support} of $\alpha$.  To prove
the claim we have two cases.  First suppose that only critical
rectangles support $\alpha$.  So $\alpha$ is normally isotopic to a
normal curve $\beta \subset \bdy \tau \cap S$.  The first step of the
tightening procedure prevents $\beta$ from being a boundary of the
almost normal piece of $S$.  It follows that $\alpha$ must be short.

Otherwise $\alpha_1$, say, is on the boundary of a critical hexagon $h
\subset f$. Let $\beta$ be a normal curve of $S$ meeting $h$ and
let $\beta_1 \subset \beta$ be one of the normal arcs in $\bdy h$. Let
$e$ be the edge of $f$ which $\alpha_1$ does not meet.  This edge is
partitioned into three pieces; $e_h \subset h$, $e'$, and $e''$. We
may assume that $\beta_1$ separates $e_h$ from $e'$.  See
Figure~\ref{Fig:AlphaShort}.

\begin{figure}
\psfrag{alpha1}{$\alpha_1$}
\psfrag{beta1}{$\beta_1$}
\psfrag{e'}{$e'$}
\psfrag{eh}{$e_h$}
\psfrag{e''}{$e''$}
\psfrag{h}{$h$}
$$\begin{array}{c}
\epsfig{file=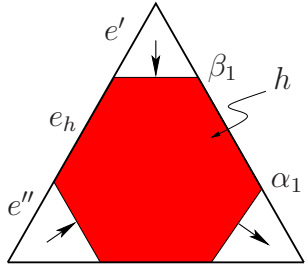, height = 3.5 cm}
\end{array}$$
\caption{The normal arcs $\alpha_1$ and $\beta_1$ are on the boundary
  of the critical hexagon $h$.  Note that $\beta$ does not meet $e'$
  or $\interior(h)$.}
\label{Fig:AlphaShort}
\end{figure}

Note that a normal curve of length $\leq 8$ has no parallel normal
arcs in a single face.  Thus $\beta$ meets $e'$ exactly once at an
endpoint of $e'$.  Since $\alpha$ and $\beta$ do not cross it follows
that $\beta$ separates $\alpha$ from $e'$ in $\bdy \tau$.

Similarly, $\alpha$ is separated from $e''$. Thus $\alpha$ does not
meet $e$ at all.  By Lemma~\ref{Lem:ShortCurves} the normal curve
$\alpha$ is short.  This finishes the proof of the claim. \qed

\begin{proofclaim}
The component $P \subset \tau \cap F_n$ has at most one boundary
component which is a normal curve.
\end{proofclaim}

Proving this will complete the lemma.  So suppose that $\bdy P$
contains two normal curves: $\alpha$ and $\beta$.  Let $A$ be the
annulus cobounded by $\alpha$ and $\beta$ in $\bdy \tau$, the boundary
of the model tetrahedron.  

Suppose now that the transverse orientation $F_n$ induces on $\alpha$
points away from $A$.  Thus $A$ and the support of $\alpha$ intersect.
There are several cases to examine, depending on the length of
$\alpha$ and the components of the support of $\alpha$.
\begin{enumerate}
\item
Suppose $\alpha$ has length three:
\begin{enumerate}
\item
If only critical rectangles support $\alpha$ then a normal triangle of
$S$ separates $\alpha$ and $\beta$.
\item
If one critical hexagon and two critical rectangles support $\alpha$
then the almost normal octagon and the exceptional tightening disk
together separate $\alpha$ and $\beta$.  (See left hand side of
Figure~\ref{Fig:WrappingThree}.)
\item
If two critical hexagons and one critical rectangle support $\alpha$
then either a normal triangle or normal quad of $S$ separates $\alpha$
and $\beta$.  (See right hand side of Figure~\ref{Fig:WrappingThree}.)
\item
If only critical hexagons support $\alpha$ then a normal triangle of
$S$ separates $\alpha$ and $\beta$.
\end{enumerate}
\item
Suppose $\alpha$ has length four:
\begin{enumerate}
\item
If only critical rectangles support $\alpha$ then a normal quad of $S$
separates $\alpha$ and $\beta$.
\item
If one critical hexagon and three critical rectangles support $\alpha$
then $S$ could not have been an almost normal surface.  (See left
hand side of Figure~\ref{Fig:WrappingFour}.)
\item
If two critical hexagons and two critical rectangles support $\alpha$
then a normal triangle of $S$ separates $\alpha$ and $\beta$.  (See
right hand side of Figure~\ref{Fig:WrappingFour}.)
\end{enumerate}
\end{enumerate}

\begin{figure}
\psfrag{alpha}{$\alpha$}
$$\begin{array}{cc}
\epsfig{file=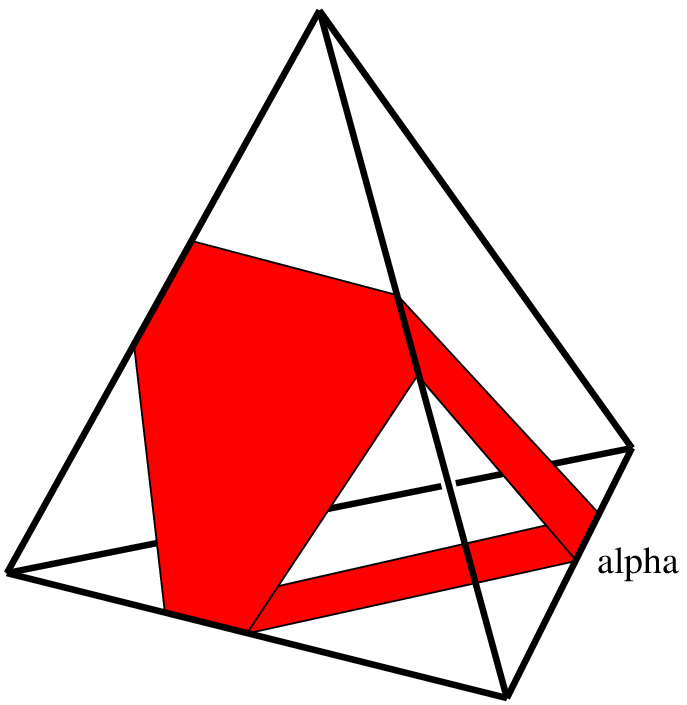, height = 3.5 cm} &
\epsfig{file=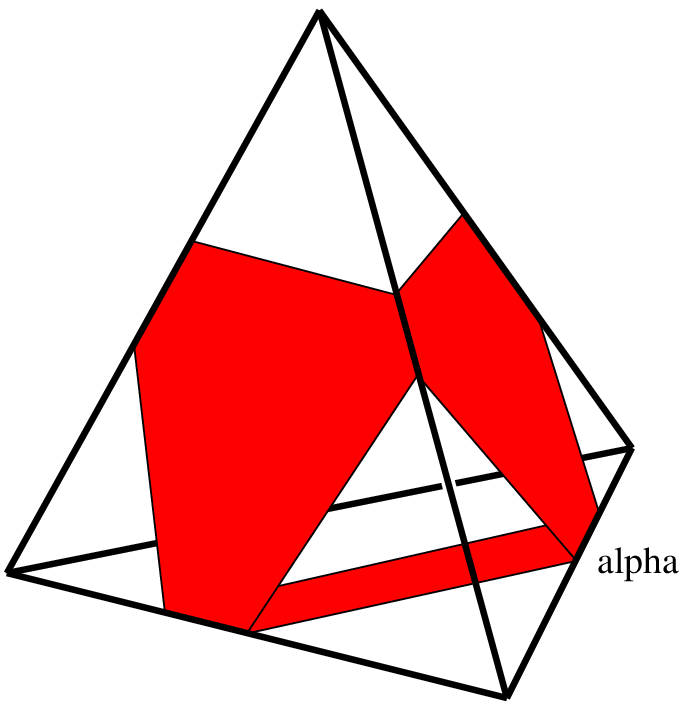, height = 3.5 cm} \\
\mbox{One hexagon} & \mbox{Two hexagons} \\
\end{array}$$
\caption{Diagrams for cases (1b) and (1c).}
\label{Fig:WrappingThree}
\end{figure}

\begin{figure}
\psfrag{alpha}{$\alpha$}
$$\begin{array}{cc}
\epsfig{file=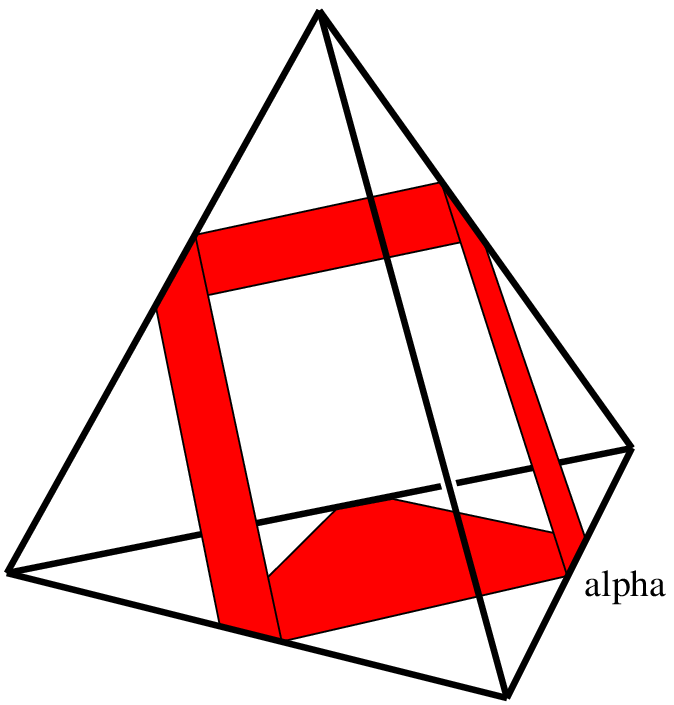, height = 3.5 cm} &
\epsfig{file=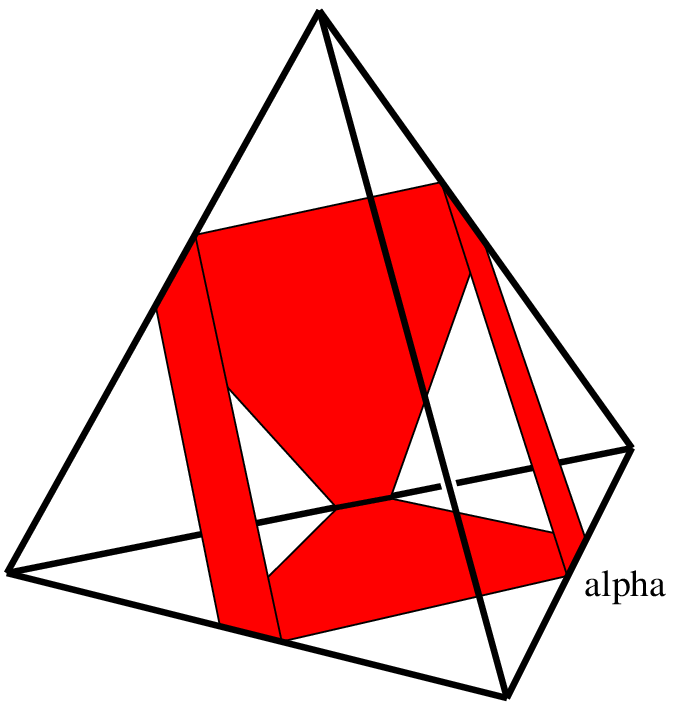, height = 3.5 cm} \\
\mbox{One hexagon} & \mbox{Two hexagons} \\
\end{array}$$
\caption{Diagrams for cases (2b) and (2c).}
\label{Fig:WrappingFour}
\end{figure}

When $\alpha$ has length four it cannot be supported by more than two
critical hexagons.

To recap: in all cases except 1(b) and 2(b), the support of $\alpha$
(possibly together with a terminal rectangle or hexagon) closes up,
implying the existence of a normal disk of $S$ with boundary a core
curve of the annulus $A$.  As this disk lies in $S$ observe that $S
\cap P \neq \emptyset$ and thus $S \cap F_n \neq \emptyset$.  This
contradicts the fact that $\calF$ is an embedding
(Lemma~\ref{Lem:Embedded}).  Case 1(b) is similar, except that the
support of $\alpha$ meets other critical or terminal components to
form the octagon piece of $S$.  So $P$ must intersect either $S$ or
the exceptional tightening disk, again a contradiction of
Lemma~\ref{Lem:Embedded}.  Lastly, in case 2(b), $S$ could not have
been almost normal.

So deduce that the transverse orientation which $F_n$ gives $\alpha$
must point toward $A$.  Thus $A$ and the support of $\alpha$ are
disjoint.  Let $\gamma$ be an arc which runs along $P$ from $\alpha$
to $\beta$.  Let $\alpha'$ be a push-off of $\alpha$ along $A$,
towards $\beta$.  This push-off bounds a disk in one of the components
of $\tau \setminus \calF$, by Lemma~\ref{Lem:Balls}. This disk does
not intersect $P \subset F_n \subset \calF$ and hence fails to
intersect $\gamma$.  This is a contradiction.
\end{proof}

\begin{remark}
\label{Rem:NestedTubes}
By Lemma~\ref{Lem:Embedded} all simple curves of $F_i$ are innermost.
It follows that the ``tubes'' analyzed in Lemma~\ref{Lem:BoundaryFn}
do not run through each other.  In addition, analysis similar to that
needed for Lemma~\ref{Lem:Handlebody} implies that these tubes are
unknotted, but this last fact is not needed in the sequel.
\end{remark}

\section{Capping off}
\label{Sec:CappingOff}

Here we construct our candidate for $C_S$, the canonical compression
body.

Let $\calF \subset |T|$ be the image of the map constructed above.
Recall that $\bdy \calF = S \cup F_n$ where $S$ is the almost normal
surface we started with and $F_n$ is the surface obtained by
``tightening'' $S$.  Note that, since $\calF$ is the embedded image of
$S \cross [0, n]$, in fact $F_n$ is isotopic to $S$ in $|T|$.  (It
cannot be normally isotopic as it has lower weight.)

\begin{define}
A two-sphere which is embedded in $|T|$ but disjoint from $T^2$ is
called a {\em bubble}.
\end{define}

We have a corollary which is easy to deduce from
Lemma~\ref{Lem:BoundaryFn}, Corollary~\ref{Cor:PlanarSurfaces}, and
Remark~\ref{Rem:NestedTubes}:

\begin{corollary}
\label{Cor:SurgeryOnFn}
Let $F'_n$ be the surface obtained by surgering all simple curves of
$F_n$.  Then $F'_n$ is a disjoint collection of bubbles and normal
surfaces.  Each bubble bounds a ball with interior disjoint from $T^2
\cap F'_n$.
\end{corollary}

Construct $C_S$ as follows: For every simple curve $\alpha$ of $F_n$
attach a two-handle to $\calF$ along $\alpha$.  Attach so that the
core of the two-handle is the subdisk of $T^2$ cut out by $\alpha$.
Call this $\calF'$.  As noted in Remark~\ref{Rem:NestedTubes} all
simple curves of $F_n$ are innermost.  So $\calF'$ is an embedded
compression body.  At this point there may be components of $\bdy_-
\calF'$ which are not normal.  By Corollary~\ref{Cor:SurgeryOnFn} all
of these are bubbles bounding a ball disjoint from all of the other
bubbles.  So cap off each bubble to obtain $C_S$.  Set $\normal{S} =
\bdy_- C_S$.  The next section proves that $v(\normal{S})$ does not
depend on the choices made in the construction of $\calF$.

\begin{remark}
\label{Rem:CompressionBodyDef}
The reason why two-spheres are allowed in $\bdy_- C_S$ should now be
clear:  we cannot prevent normal two-spheres from appearing in the
normalization process.  
In particular, if $S$ is an almost normal two-sphere then, for one of
the two possible transverse orientations, there will always be a
normal two-sphere appearing in $\normal{S}$.
\end{remark}

\section{Proof of the normalization theorem}
\label{Sec:NormalizationAlg}

Suppose that $S$ is almost normal and equipped with a transverse
orientation.  Before proving Theorem~\ref{Thm:NormalizationAlg} recall
that $C_S$, a compression body in $|T|$, is {\em canonical} for $S$ if
$\bdy_+ C_S = S$, $\bdy_- C_S$ is normal, the transverse orientation
on $S$ points into $C_S$, and any normal surface $S' \subset |T|$ may
be normally isotoped to one disjoint from $C_S$.

We now have:

\begin{theorem}
\label{Thm:NormalizationAlg}
Given a transversely oriented almost normal surface $S$ there exists a
canonical compression body $C_S$ and $C_S$ is unique (up to normal
isotopy).  Furthermore there is a algorithm which, given the
triangulation $T$ and the surface vector $v(S)$, computes the surface
vector of $\bdy_- C_S = \normal{S}$.
\end{theorem}

\begin{proof}

We proceed in several steps.

\begin{proofclaim}
A canonical compression body $C_S$ exists.
\end{proofclaim}

There are two cases.  Either the transverse orientation for $S$ points
at the exceptional surgery disk (implying that $S$ contained an almost
normal annulus) or the transverse orientation points at an exceptional
tightening disk.

In the first case, $C_S$ is obtained by thickening $S$ slightly and
adding a regular neighborhood of the exceptional surgery disk.  It is
clear that $C_S$ is a compression body, $\bdy_+ C_S = S$, and $\bdy_-
C_S$ is normal.  Suppose that $S'$ is any normal surface in $T$ which
is disjoint from $S$.  Then, perhaps after a normal isotopy of $S'$
(rel $S$), we have that $S'$ is disjoint from the exceptional surgery
disk for $S$.  It follows that $S'$ may be isotoped out of $C_S$.

In the second case the transverse orientation of $S$ points at an
exceptional tightening disk of $S$.  As in Section~\ref{Sec:Tight}
form $\calF$ with $\bdy \calF = S \cup F_n$.  As in
Section~\ref{Sec:CappingOff} attach two-handles to $\calF$ along the
simple curves of $F_n$ to obtain $\calF'$.  Cap off the bubbles with
their three-balls to obtain $C_S$.  Again, $C_S$ is a compression body
with $\bdy_+ C_S = S$.

Suppose now that $S'$ is some normal surface in $T$ which is disjoint
from $S$.  Then, by Corollary~\ref{Cor:Barrier}, the surface $S'$ is
disjoint from $\calF$ (perhaps after a normal isotopy of $S'$ rel
$S$).  Since $S'$ is normal it cannot meet any of the disks (in $T^2$)
bounded by simple curves of $F_n$.  So $S' \cap \calF' = \emptyset$ as
well.  Finally, suppose that $A$ is a bubble component of $\bdy_-
\calF'$.  Let $B$ be the three-ball which $A$ bounds (such that $B
\cap T^2 = \emptyset$).  Then no component of $S'$ meets $B$ as $S'
\cap A = \emptyset$ and $S'$ is normal.  Deduce that $S' \cap C_S =
\emptyset$.  The claim is complete.

\begin{proofclaim}
The canonical compression body $C_S$ is unique (up to normal
isotopy). 
\end{proofclaim}

Suppose that $C_S$ and $C'_S$ are both canonical compression bodies
for $S$.  Let $A = \bdy_- C_S$ and $A' = \bdy_- C'_S$.  Then $A$ and
$A'$ are normal surfaces, both disjoint from $S$.  It follows that
there exists a normal isotopy $\calH$ which moves $A'$ out of $C_S$,
rel $S$, and conversely another normal isotopy $\calH'$ which moves
$A$ out of $C'_S$, rel $S$.

Consider any face $f \in T^2$ and any normal arc $\alpha \subset f
\cap S$.  Let $X \subset f \cap C_S$ be the component containing
$\alpha$.  Also take $X'$ to be the component of $f \cap C'_S$ which
contains $\alpha$.  We must show that $X$ and $X'$ have the same
combinatorial type.  Suppose not.  After possibly interchanging $X$
and $X'$ there are only six situations to consider:

\begin{enumerate}
\item $X$ is a critical rectangle and $X'$ is a terminal rectangle.
\item $X$ is a critical rectangle and $X'$ is a critical hexagon.
\item $X$ is a critical rectangle and $X'$ is a terminal hexagon.
\item $X$ is a critical hexagon and $X'$ is a terminal hexagon.

In any of these four cases let $\delta$ be the normal arc of $A =
\bdy_- C_S$ on the boundary of $X$.  Note that $\bdy X'$ contains
$\alpha$ (as does $\bdy X$) and also another normal arc $\beta \subset
f \cap S$ which does not meet $X$ (as $S = \bdy_+ C_S$).  Now note
that it is impossible for $\calH'$ to normally isotope $\delta$ out of
$X'$ while keeping $S$ fixed pointwise (as $\delta$ would have to
cross $\beta$).

\item $X$ is a terminal rectangle and $X'$ is a critical hexagon.
\item $X$ is a terminal rectangle and $X'$ is a terminal hexagon.

In either of these cases let $\beta$ be the other normal arc of $S
\cap \bdy X$.  Then $\beta$ intersects the interior of $X'$, a
contradiction.
\end{enumerate}

This proves the claim.

\begin{proofclaim}
There is a algorithm which, given the triangulation $T$ and the
surface vector $v(S)$, computes the surface vector of $\bdy_- C_S =
\normal{S}$.
\end{proofclaim}

We follow the proof of Lemma~\ref{Lem:Embedded}: We keep track of the
intersection between the image of $\calF_k$ and every face $f \in
T^2$.  These are unions of components, with all allowable kinds shown
(up to symmetry) in Figures~\ref{Fig:Rectangles}
and~\ref{Fig:Hexagons}.  There is at most one hexagon in each face and
perhaps many rectangles, arranged in three families, one for each
vertex of $f$.  At stage $n$ there are no bent arcs remaining.  Now
delete all simple curves of $F_n$ and all normal arcs of $S$.  The
normal arcs left completely determine $\normal{S}$ and from this we
may find the surface vector $v(\normal{S})$.  This proves the claim
and finishes the proof of Theorem~\ref{Thm:NormalizationAlg}.
\end{proof}

Of course, the algorithm just given is inefficient.  It depends
polynomially on $\size(T)$ and $\weight(S)$.  In the next section we
improve this to a algorithm which only depends polynomially on
$\size(T)$ and $\log(\weight(S))$.

As a corollary of Theorem~\ref{Thm:NormalizationAlg}:

\begin{corollary}
\label{Cor:NormalizingS2}
If $S \subset |T|$ is a transversely oriented almost normal two-sphere
then $C_S$ is a three-ball, possibly with some open three-balls
removed from its interior.  (These have closures disjoint from each
other and from $S$.)
\end{corollary}

Now an orientable surface in an orientable three-manifold may be
transversely oriented in exactly two ways.  By
Theorem~\ref{Thm:NormalizationAlg}, if $S$ is an almost normal
surface, for each transverse orientation there is a canonical
compression body.  Call these $C_S^+$ and $C_S^-$.

From Corollary~\ref{Cor:NormalizingS2} deduce:

\begin{theorem}
\label{Thm:NormalizingS2}
If $S \subset |T|$ is an almost normal two-sphere and both $\bdy C_S^+
\setminus S$ and $\bdy C_S^- \setminus S$ are (possibly empty)
collections of vertex-linking two-spheres, then $|T|$ is the
three-sphere.
\end{theorem}

\begin{proof}
By hypothesis $\bdy C_S^+ \setminus S$ is a collection of vertex
linking spheres.  For each of these add to $C_S^+$ the corresponding
vertex neighborhood.  Let $\BB^+$ be the resulting submanifold of
$|T|$.  By Alexander's Trick $\BB^+$ is a three-ball.  Do the same to
$C_S^-$ to produce $\BB^-$.  Applying Alexander's Trick again deduce
that the manifold $|T| = \BB^+ \cup_S \BB^-$ is the three-sphere.
\end{proof}

\section{An example}
\label{Sec:Example}

\begin{figure}
\psfrag{0}{$0$}
\psfrag{1}{$1$}
\psfrag{2}{$2$}
\psfrag{3}{$3$}
$$\begin{array}{c}
\epsfig{file=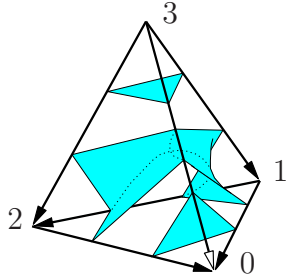, height = 3.5 cm}
\end{array}$$
\caption{A one tetrahedron triangulation of $S^3$.  It is a simple
  exercise to list all normal and almost normal surfaces in $T$.  The
  reader may further amuse herself by drawing the graph $T^1$ as it
  actually sits in $S^3$.  Slightly harder is to draw the
  two-skeleton.}
\label{Fig:Example}
\end{figure}

Here we give a brief example of the normalization procedure.  Let $T$
be the one vertex triangulation shown in Figure~\ref{Fig:Example}.

The front two faces ($1$ and $2$) are glued to each other as are the
back faces ($0$ and $3$).  The faces are glued to give the edge
identifications shown.  The surface $S$ depicted in $T$ is an almost
normal two-sphere with two triangles and one almost normal octagon.
It is easy to check this by computing $\euler(S) = 3 - 7 + 6 = 2$.

The sphere $S$ has two exceptional tightening disks: $D$ meeting the
edge $(0,3)$ of the model tetrahedron and $D'$ meeting edge $(1,2)$.

Tightening along $D$ gives $F_1$ which is the vertex link.  Tightening
along $D'$ gives $F_1', F_2', F_3'$.  Here $F_3'$ is a weightless
two-sphere in $T$ with a single simple curve and no other intersection
with the two-skeleton.  As a note of caution: $F_1'$ drawn in the
model tetrahedron has four bent arcs -- however $F_1' \cap T^2$
contains only two.  These are independent of each other and doing
these moves in some order gives $F_2'$ and $F_3'$.  To obtain the
normalization of $S$ on the $D'$ side, surger the simple curve of
$F_3'$ and cap off the two resulting bubbles.  

So, on the $D$ side of $S$ the normalization is the vertex link.  On
the $D'$ side the normalization is the empty set.  It follows from
Theorem~\ref{Thm:NormalizingS2} that $|T|$ is the three-sphere.
This finishes the example.

\section{Normalizing quickly}
\label{Sec:NewNormalization}

The normalization procedure can be made much more efficient. 

\begin{theorem}
\label{Thm:NewNormalizationAlg}
There is a polynomial time algorithm which, given $T$ and the surface
vector $v(S)$, produces as output $v(\normal{S})$, the normalization
of $S$.  Here $S$ is assumed to be a separating transversely oriented
almost normal surface and $T$ is assumed to be a triangulation
of a three-manifold.
\end{theorem}

\begin{remark}
\label{Rem:NewNormalizationAlgNonseparating}
In fact $S$ need not be separating.  However the proof is somewhat
simpler in the separating case and is all we require in what follows.
\end{remark}

Recall that $N_S$ is the closure of the component of $|T| \setminus S$
which the transverse orientation points into.  Then $\blocky{N}_P$ is
the union of all product blocks in $N_S$ and $\blocky{N}_C$ is the
union of all the core blocks.  Also $N_P$ is a regular neighborhood of
$\blocky{N}_P$, taken in $N_S$.  Finally $N_C = \closure{N_S \setminus
N_P}$.  We will prove Theorem~\ref{Thm:NewNormalizationAlg} by
altering our original normalization procedure three times.  First we
will show that the order of the tightening moves is irrelevant.  Then
we show that surgeries on simple curves and capping off of bubbles may
happen during the normalization process, instead of being held until
the end.  Finally we show that tightening inside of $N_P$ can be done
very quickly.  These three modifications combine to give an efficient
algorithm.

\subsection{Changing the order of the tightening moves}

As stated in Remark~\ref{Rem:CalFNotUnique} the isotopy $\calF \from S
\cross [0,n] \to M$ need not be unique.  But we do have:

\begin{lemma}
\label{Lem:MoveOrder}
Any order for the tightening moves (performed in the construction of
$\calF$) gives the same surface $\normal{S}$ once the simple curves of
$F_n$ have been surgered.
\end{lemma}

This follows immediately from the first sentence of
Theorem~\ref{Thm:NormalizationAlg}.

\subsection{Surgery on simple curves and deleting bubbles}

We now alter the tightening procedure in a more substantial fashion:

Recall that $S \subset |T|$ is a transversely orientable separating
almost normal surface.  Recall that $D$ is the exceptional tightening
disk for $S$.  Transversely orient $S$ to point into the component of
$|T| \setminus S$ which meets $D$.  Here is the {\em $G$-tightening
procedure}:

\begin{enumerate}
\item
Let $G_0 = S$.  Let $D_0 = D$.

\item
Now do a small normal isotopy of $G_0$ in the transverse direction
while also tightening $G_0$ along $D_0$.  Call the surface so obtained
$G'_0$.  Now surger all simple curves of $f \cap G'_0$ for every $f
\subset T^2$ to obtain $G''_0$.  Then delete any bubble components of
$G''_0$ (\ie, two-sphere components which are contained in the
interior of tetrahedra).  Call the resulting surface $G_1$.  Note that
$G_1$ inherits a transverse orientation from $G_0$.

\item
At stage $k \geq 1$ there are two possibilities.  Suppose first that
$G_k$ has an outermost bent arc $\alpha$ with the transverse
orientation of $G_k$ pointing into the tightening disk $D_k$, which is
cut out of $T^2$ by $\alpha$.  Then perform a small normal isotopy of
$G_k$ in the transverse direction while tightening $G_k$ across $D_k$.
Call the surface so obtained $G'_k$.  Now surger all simple curves of
$f \cap G'_k$ for every $f \in T^2$ to obtain $G''_k$.  Then
delete any bubble components of $G''_k$.  Call the resulting surface
$G_{k+1}$.  Note that $G_{k+1}$ inherits a transverse orientation from
$G_k$.

If there is no such outermost bent arc $\alpha \subset G_k$ then set
$n = k$ and the procedure halts.
\end{enumerate}

\begin{lemma}
\label{Lem:SurgeryOrder}
The surface $G_n$ is normally isotopic to $\normal{S}$, the
normalization of $S$.
\end{lemma}

\begin{proof}
Recall that Lemma~\ref{Lem:Embedded} gives a complete classification
of the possible components of intersection of $\image(\calF_k)$ with
the faces of $T^2$.  Again, see Figures~\ref{Fig:Rectangles}
and~\ref{Fig:Hexagons}.  The only components containing a simple curve
are the terminal rectangle with hole and terminal hexagon with hole.
Hence their names.

Since the terminal with hole rectangles and hexagons do not contain
normal or bent arcs of $F_k$ they remain unchanged in the
$F$-tightening procedure until $F_n$ is reached.  Then all simple
curves are surgered and bubbles capped off.  Thus it makes no
difference to the resulting surface $\normal{S}$ if we surger simple
curves and delete bubbles as soon as they appear.
\end{proof}

\subsection{Tightening in I-bundle regions}

We now give the final modification of the tightening procedure.
Suppose that $v(S)$ is an almost normal surface vector.  Suppose also
that $S$ has a transverse orientation pointing at an exceptional
tightening disk.

Recall that $N_S$ is the blocked submanifold cut from $|T|$ by the
surface $S$ (so that the transverse orientation points into $N_S$).
Also, $N_P$ is the $I$-bundle part of $N_S$ while $N_C = \closure{N_S
\setminus N_P}$ is the core of $N_S$.

We require slightly more sophisticated data structures.  First
define $\product(S)$ to be the list $\{v_j\}_{j=1}^m$ where the
$j^\th$ element is the vector $2 \cdot v(N_P^j)$ -- here $v(N_P^j$) is
the block vector for the $j^\th$ component of $\blocky{N}_P$, found by
Theorem~\ref{Thm:AHTBlockVectorAlg}.  That is, $\sum v_j$ counts the
normal disks of $S$ which make up the horizontal boundary of the
product blocks in $N_S$.

Put a copy of the horizontal boundary of $N_C$ in $\core(S)$.  That
is, record in $\core(S)$ all of the gluing information between edges
of disks which are on the horizontal boundary of core blocks.  Also
record, for each edge in $\bdy \vbdy N_C$, which disk of $\core(S)$
contains it and which component $N_P^j \subset N_P$ it is glued to.
Build a model of $N_C$. That is, deduce what core blocks occur in
which tetrahedra and how they are glued across faces of $T^2$.

\begin{figure}
\psfrag{D}{$D_k$}
\psfrag{NP}{$N_P^j$}
\psfrag{bdy}{$\vbdy N_P^j$}
$$\begin{array}{c}
\epsfig{file=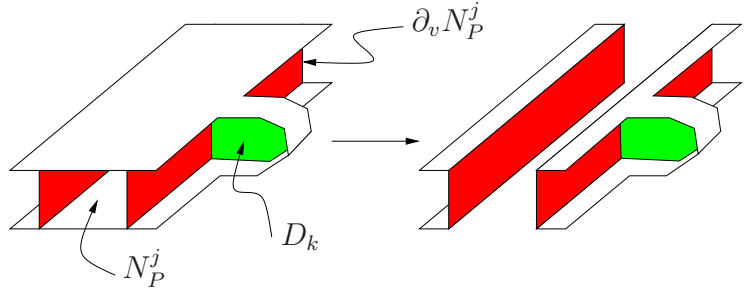, height = 3.5 cm}
\end{array}$$
\caption{Removing the horizontal boundary of $N_P^j$ and adding the vertical.}
\label{Fig:ProductSurgery}
\end{figure}

We now turn to constructing a sequence of surfaces $H_k$.  Each $H_k$
will be represented by the two pieces of data: $\core(H_k)$ and
$\product(H_k)$.  Here is the {\em $H$-tightening procedure}:

\begin{enumerate}
\item
Let $\core(H_0) = \core(S)$.  Let $\product(H_0) = \product(S)$.  Let
$D_0 = D$.  At stage $k$ there is a tightening disk $D_k$ used to
alter $H_k$.

\item
If the $D_k$ has empty intersection with $N_P$ then perform the
tightening move as in the $G$-sequence.  This effects only the pieces
in $\core(H_k)$ and we use the tightening move to compute
$\core(H_{k+1})$.  Set $\product(H_{k+1}) = \product(H_k)$ and go to
stage $k+1$.

\item 
Suppose $D_k$ intersects a component of $N_P$, say $N_P^j$.  Then set
$\product(H_{k+1}) = \product(H_K) \setminus \{v_j\}$; \ie, remove
$v_j$ from the product part.  We also alter the disks in the core as
follows: Let $\core(H'_k) = \core(H_k) \cup \vbdy N_P^j$.  Let $D'_k =
D_k \setminus N_P^j$ (that is, remove a small neighborhood of $T^1$
from $D_k$).  See Figure~\ref{Fig:ProductSurgery}.  Then $D'_k$ is a
surgery disk for $\core(H'_k)$.  So surger along $D'_k$, surger along
all simple curves of $\core(H'_k)$, and delete all bubbles in
$\core(H'_k)$.  This finally yields $\core(H_{k+1})$. Go on to stage
$k+1$.

\item
At stage $k+1$ there are two possibilities: either there is a bent arc
in $\core(H_{k+1})$ or there is not.  If there is then we have a
tightening disk $D_{k+1}$ and proceed as above.  If there is no bent
arc then sum the vectors in $\product(H_{k+1})$ and add to this vector
the number of normal disks of each type in $\core(H_{k+1})$.  Output
the final sum $v(H_n)$.
\end{enumerate}

This is our final modification of the tightening procedure.

\subsection{Correctness and efficiency}

\begin{proof}[Proof of Theorem~\ref{Thm:NewNormalizationAlg}]

Note that if the transverse orientation on $S$ points towards an
exceptional surgery disk of $S$ then the theorem is trivial.  So
suppose instead that a tightening disk is pointed at.

\begin{proofclaim}
The $H$-tightening procedure outputs $v(\normal{S})$.
\end{proofclaim}

That is, we claim that $H_n$ is normally isotopic to $\normal{S}$.
Since the $H$ procedure is identical to the $G$ procedure in $N_C$ we
need only consider the situation where a tightening disk meets $N_P$.
Consider the smallest $k$ such that $D_k \cap N_P \neq \emptyset$.
Recall that $\bdy D_k = \alpha \cup \beta$ where $\alpha \subset H_k$
and $\beta \subset T^1$. Also by our hypothesis on $k$ the arc $\beta$
is contained in $T^1 \cap \vbdy \blocky{N}_P$ while only a small
neighborhood of $\bdy \alpha$ (taken in $\alpha$) is contained in
$N_P$.  Suppose that $N_P^j$ is the component of $N_P$ containing
$\beta$.

We have assumed that $H_k = G_k$.  We will show that we can reorder
tightening moves in the $G$-procedure to obtain $G_{k+k'}$
normally isotopic to $H_{k+1}$.  Then it will follow from
Lemma~\ref{Lem:SurgeryOrder} that the $H$ procedure produces correct
output.

Recall that $\blocky{N}_P^j$ and $N_P^j$ are $I$-bundles.  Let $\pi$
be the natural quotient map which squashes $I$-fibres to a point.  Let
$E = \pi(N_P^j)$.  Let $\blocky{E} = \pi(\blocky{N}_P^j)$.  Note that
$\blocky{E}$ is not necessarily a surface.  However $E$ is a surface
with boundary, $\blocky{E}$ naturally embeds in $E$, and there is a
small deformation retraction of $E$ to $\blocky{E}$.  Note that
$\blocky{E}$ and $E$ inherit cell structures from $\blocky{N}_P^j$ and
$N_P^j$.  Choose a spanning tree $U$ for the one-skeleton
$\blocky{E}^1$ of $\blocky{E}$ rooted at $b = \pi(\beta)$.  Choose an
ordering of the vertices of $U$, $\sigma \from U^0 \to (\NN \cap [1,
k'])$, so that for any vertex $d$ with parent $c$ we have $\sigma(c) <
\sigma(d)$.  Here $k' = |U^0|$ is the number of vertices in $U^0$.

We now have a sequence of tightening moves to perform in the $G$
procedure.  At step one do the tightening move along the disk $D_k$,
surger all simple curves, and delete bubbles.  At step $i, i > 1$,
examine the edge $e$ between $c$ and $d$ (where $\sigma(d) = i$ and
$c$ is the parent of $d$).  Then, by induction and the fact that
$\sigma(c) < \sigma(d) = i$ there is a bent arc of $G_{k+i-1}$ in the
rectangle $\pi^{-1}(e)$ with endpoints on the segment $\pi^{-1}(d)
\subset T^1$.  Do this tightening move, surger simple curves, delete
bubbles, and go to stage $i+1$.

After stage $k' = |U^0|$ we obtain the surface $G_{k + k'}$ which is
normally isotopic to the following: $(G_k \setminus \hbdy N_P^j) \cup
\vbdy N_P^j$ surgered along the disk $D'_k$, surgered along simple
curves, with bubbles deleted.  Here $D'_k = \closure{D_k \setminus
N_P^j}$.  So $G_{k+k'}$ agrees with $H_{k+1}$ and the claim is proved.

\begin{proofclaim}
Precomputation for the $H$ procedure takes time at most polynomial in
$\size(T)$ and $\log(\weight(S))$.
\end{proofclaim}

Theorem~\ref{Thm:AHTBlockVectorAlg} computes the vectors
$\{v_j\}_{j=1}^m$ in the required amount of time.  This gives
$\product(N_S)$.  Then, since there are only a linear number (in
$\size(T)$) of core blocks in $N_S$ (Remark~\ref{Rem:BoundedCore}) we
can also compute their gluings and so compute $\core(S)$ in the
alloted time.

\begin{proofclaim}
The number of steps in the modified normalization procedure is
polynomial in $\size(T)$.
\end{proofclaim}

Each step either reduces the weight of $\core(H_k)$ by two or removes
a vector from $\product(H_k)$.  Since the weight of $\core(H_k)$ is at
most linear (again, Remark~\ref{Rem:BoundedCore}), and since there are at
most a linear number of components of $N_P$ (see
Remark~\ref{Rem:NumberOfN_PComponents}), the claim follows.

\begin{proofclaim}
Performing each step of the modified normalization procedure takes
time at most polynomial in $\size(T)$ and $\log(\weight(S))$.
\end{proofclaim}

If the tightening disk is disjoint from $N_P$ then we only have to
alter $\core(H_k)$ in the tetrahedra adjacent to the disk.  There are
only a linear number of these.

If the tightening disk meets a component of $N_P$, say $N_P^j$, then
delete $v_j$ from $\product(H_k)$ in polynomial time (in
$\size(T)$).  Alter $\core(H_k)$ by gluing on a copy of
$\vbdy N_P^j$, surgering along the remnants of the tightening disk
$D_k$, surgering all simple curves, and deleting bubbles.  As $\vbdy
N_P^j$ is a subset of $\vbdy N_C$ it is at most linear in size (in
terms of $\size(T)$).  Thus we can make the desired changes in the
required time.

To sum up: we can compute the desired result,
$v(\normal{S})$, in time which is at most a product of
polynomials in $\size(T)$ and $\log(\weight(S))$.  This completes the
proof of Theorem~\ref{Thm:NewNormalizationAlg}.
\end{proof}

\section{Crushing, or: ``New triangulations for old''}
\label{Sec:Crushing}

Let $T$ be a triangulation of a closed three-manifold.
Suppose we are given a choice of quad type in a single tetrahedron,
say the $a^\th$ type of quad in $\tau_i$.  Here $a \in \{1, 2, 3\}$
and the other two elements of $\{1, 2, 3\}$ are $b$ and $c$.  Recall
that the $a^\th$ quad type separates the vertices $0$ and $a$ from the
vertices $b$ and $c$.

Let $\theta$ be the permutation $(0a)(bc)$.  Let $\{(i, j_s,
\sigma_s)\}_{s=0}^3$ be the four face pairings with $i$ as the first
element.  Here $\sigma_s$ glues the $s^\th$ face of $\tau_i$ to some
face of $\tau_{j_s}$.  Note that $\{(j_s, i, \sigma_s^{-1})\}_{s=0}^3$
are also face pairings in $T$.

Define a new triangulation $T'$ by {\em crushing} the tetrahedron
$\tau_i$ along the $a^\th$ quad, as follows: Delete $\tau_i$ from $T$.
Delete all of the face pairings $\{(i, j_s, \sigma_s)\}_{s=0}^3$.
Replace the face pairing $(j_s , i, \sigma_s^{-1})$ (if $i \neq j_s$)
with 
$$\left( j_s, j_{\theta(s)}, 
       \sigma_{\theta(s)} \cdot R_{(s, \theta(s))} \cdot
       \sigma_s^{-1} \right),$$
for $s \in \{0, 1, 2, 3\}$.  Here $R_{(s, \theta(s))}$ is the rotation
of the model tetrahedron, about the edge with vertices $\{0, 1, 2, 3\}
\setminus \{s, \theta(s)\}$, which takes face $s$ to face $\theta(s)$.
Note that $\sigma_{\theta(s)} \cdot R_{(s, \theta(s))} \cdot
\sigma_s^{-1}$ is the composition of three orientation reversing maps
and thus is also orientation reversing.  Finally, no face of any model
tetrahedron in $T'$ is glued to itself -- thus $T'$ is a
triangulation.

To keep track of this operation it may help to refer to the picture of
a quad of type $3$ shown on the right hand side of
Figure~\ref{Fig:NormalDisks}.

Now suppose that $p$ is a {\em polarization} of the triangulation $T$;
that is, $p$ is a map from the set of tetrahedra to the set $\{0, 1,
2, 3\}$.  Produce a new triangulation $T'$ by {\em crushing} $T$ along
$p$: To begin with let $T'$ be an exact copy of $T$.  Now, for each $i
= 1, 2, \ldots, \size(T)$ do one of two things; If $p(\tau_i) = 0$
simply go on to $i+1$.  If $p(\tau_i) \neq 0$ then remove $\tau_i$ by
crushing along the $p(\tau_i)$ quad, as above, and go on to $i+1$.

The following theorem is now clear:

\begin{theorem}
\label{Thm:CrushingAlg}
There is an polynomial time algorithm which, given a triangulation $T$
and a polarization $p$, produces $T'$, the triangulation of $T$
crushed along $p$. \qed
\end{theorem}

Of more interest is:

\begin{theorem}
\label{Thm:CrushingAndConnectSum}
Suppose $T$ is a triangulation so that the connect sum $\connect |T|$
is a homology three-sphere.  Suppose $p$ is a polarization
coming from $S$, a non-vertex linking normal two-sphere.  Then the
triangulation $T'$, the result of crushing $T$ along $p$, satisfies
$\connect |T'| \homeo \connect |T|$.
\end{theorem}

\begin{proof}
Theorem~4.10 of Jaco and Rubinstein's paper~\cite{JacoRubinstein02a}
essentially claims this result for any closed, orientable
three-manifold $|T|$ with the caveat that some connect summands of
$|T|$ homeomorphic to lens spaces may by omitted from the crushed
$|T'|$.  (See also Theorem~3.1 of~\cite{Barchechat03}.)

However, by Lemma~\ref{Lem:ConsequencesOfHomologySphere} no
non-trivial lens space appears as a connect summand of the
homology three-sphere $|T|$.  Finally, omitting $S^3$
summands does not change the connect sum.  The result follows.
\end{proof}

\section{Thompson's Theorem}
\label{Sec:Thompson}

We will need to use Casson's version~\cite{Casson97} of the proof of
Thompson's Theorem~\cite{Thompson94} (which in turn relies heavily on
work of Rubinstein~\cite{Rubinstein97}):

\begin{theorem}
\label{Thm:CassonAlg}
There is an exponential time algorithm which, given a triangulation $T$,
decides whether or not $|T|$ is homeomorphic to the three-sphere.
\end{theorem}

We now sketch Casson's version of Thompson's algorithm.  Begin with a
triangulation $T_0 = T$.  Check, using
Theorems~\ref{Thm:ThreeManifoldAlg} and~\ref{Thm:HomologySphereAlg},
that $T_0$ is a homology three-sphere.  Inductively we have a
triangulation $T_i$.

If $T_i$ is not zero-efficient then apply
Lemma~\ref{Lem:FundamentalSphereAlg} to find $S_i \subset |T_i|$, a
fundamental non-vertex-linking normal two-sphere.  Let $T_{i + 1}$ be
the triangulation obtained by crushing along $S_i$.  This requires
Theorem~\ref{Thm:CrushingAlg}.

If $T_i$ is zero-efficient use Lemma~\ref{Lem:FundamentalSphereAlg} to
search for almost normal two-spheres.  If some component of $T_i$ does
not contain an almost normal two-sphere then by
Theorem~\ref{Thm:CrushingAndConnectSum} and
Theorem~\ref{Thm:ThompsonSphere} the manifold $|T|$ was not the
three-sphere.  If $S_i$ is an almost normal two-sphere inside a
component $T'$ of $T_i$ then let $T_{i + 1} = T_i \setminus T'$.

This completes the algorithm.  If $T_n$ is non-empty, then $|T|$ was
not the three-sphere.  If $T_n$ is empty then $|T|$ was homeomorphic
to the three-sphere.  Both of these again use
Theorem~\ref{Thm:CrushingAndConnectSum}.  This completes our
description of Casson's algorithm and our proof of its correctness.

Note that $\size(T_i) + i \leq \size(T)$ as either crushing along a
polarization or deleting a component always reduces the number of
tetrahedra by at least one.  This finishes the sketch of the proof of
Theorem~\ref{Thm:CassonAlg}.

\section{Showing the problem lies in NP}
\label{Sec:NP}

We are now in a position to prove:

\begin{theorem}
\label{Thm:RecognitionInNP}
The three-sphere recognition problem lies in the complexity class \NP.
\end{theorem}

\begin{proof}
Suppose that $T$ is a triangulation of the three-sphere.  The
certificate is a sequence of pairs $(T_i, v(S_i))$ with the following
properties.
\begin{itemize}
\item
$T = T_0$.
\item
$S_i$ is a normal or almost normal two-sphere, contained in $|T_i|$,
with $\weight(S_i) \leq \exp(\size(T_i))$.
\item
If $S_i$ is normal then $S_i$ is not vertex linking and $T_{i+1}$ is
obtained from $T_i$ by crushing along $S_i$.
\item
if $S_i$ is almost normal then $S_i$ normalizes to vertex linking
two-spheres, in both directions.  Also, $T_{i+1}$ is obtained from
$T_i$ by deleting the component $T'$ of $T_i$ which contains $S_i$.
\item
Finally, the last triangulation $T_n$ is empty, as is $S_n$.
\end{itemize}

Note that existence of the certificate is given by our proof of
Theorem~\ref{Thm:CassonAlg}.  So the only task remaining is to check
the certificate.  Here we find two subtle points -- we will not
attempt to verify that the $S_i$ are fundamental nor will we try to
check that the $T_i$ containing almost normal two-spheres are
zero-efficient.

Note instead, since the $S_i$ are fundamental, they obey the weight
bounds given in Lemma~\ref{Lem:FundamentalBound}; that is,
$\weight(S_i) \leq \exp(\size(T_i))$.

So suppose a certificate $(T_i, v(S_i))$ as above, for the triangulation
$T$, is given to us.  First check, using
Theorem~\ref{Thm:ThreeManifoldAlg} and~\ref{Thm:HomologySphereAlg},
that $T$ is a triangulation of a homology three-sphere.

By Theorem~\ref{Thm:IdenticalAlg} check that $T = T_0$.  Using
Theorem~\ref{Thm:ConnectedAlg} verify that $S_i$ is a connected normal
or almost normal surface.  Using Lemma~\ref{Lem:EulerAlg}
compute the Euler characteristic of $S_i$.  (Here we are using the
fact that $\weight(S_i) \leq \exp(\size(T_i))$ in order to compute
Euler characteristic in time only polynomial in $\size(T_i)$.)  This
verifies that $S_i$ is a two-sphere.

If $S_i$ is normal, by Theorem~\ref{Thm:CrushingAlg}, crush $T_i$
along $S_i$ in time at most polynomial in $\size(T_i)$.  Then
check, using Theorem~\ref{Thm:IdenticalAlg}, that $T_{i + 1}$ agrees
with the triangulation obtained by crushing $T_i$.

If $S_i$ is almost normal, we need to check that $T'$, the component
of $T_i$ containing $S_i$, has $|T'| \homeo S^3$.  Using
Theorem~\ref{Thm:NewNormalizationAlg} normalize $S_i$ in both
directions in time at most polynomial in $\size(T_i)$ (again, because
$\log(\weight(S_i)) \leq \log(\exp(\size(T_i))) = \size(T_i)$).  If
all components of the two normalizations $\normal{S_i^+}$ and
$\normal{S_i^-}$ are vertex linking two-spheres then $T'$ is a
triangulation of the three-sphere, by Theorem~\ref{Thm:NormalizingS2}.
Finally, use Theorem~\ref{Thm:IdenticalAlg} to check that the
triangulation $T_i \setminus T'$ is identical to $T_{i+1}$.
\end{proof}

\section{Questions and future work}
\label{Sec:Questions}

Our techniques should also apply to the following question:

\begin{question}
\label{Que:SurfaceBundles}
Is the surface bundle recognition problem in \NP?
\end{question}

Given a triangulation $T$ of a surface bundle the certificate would be
a certain normal two-sphere $S_0$ and collections of surfaces
$\{F_i\}$ and $\{G_i\}$.  Here all of the $F_i$ and $G_i$ are fibres
of some surface bundle structure on $|T|$, the $F_i$ are all normal,
the $G_i$ are all almost normal, $G_i$ normalizes to $F_i$ and
$F_{i+1}$ for $i > 0$, and $G_0$ normalizes to $F_0 \cup S_0$ on one
side and to $F_1$ on the other.  (See~\cite{Schleimer01} for some of
the necessary existence results.)  Also, $S_0$ bounds a three-ball
$B_0$ in $|T|$ which contains all of the vertices of $T$.  (Note that
there is a non-trivial issue here: Corollary~\ref{Cor:ThreeBallInNP}
must be modified to allow us to certify that $B_0$ is a three-ball.)


Given this it is not unreasonable to ask:

\begin{question}
\label{Que:SurfaceBundleNPHard}
Is the surface bundle recognition problem \NP-hard?
\end{question}

Perhaps more difficult to resolve would be:

\begin{question}
\label{Que:ThreeSphereNPHard}
Is the three-sphere recognition problem \NP-hard?
\end{question}

\bibliographystyle{plain}
\bibliography{bibfile}

\begin{thebibliography}{10}

\bibitem{AgolEtAl02}
Ian Agol, Joel Hass, and William~P. Thurston.
\newblock The computational complexity of knot genus and spanning area.
\newblock
  \href{http://front.math.ucdavis.edu/math.GT/0205057}{arXiv:math.GT/0205057}.

\bibitem{Barchechat03}
Alexandre Barchechat.
\newblock {\em Minimal triangulations of 3-manifolds}.
\newblock PhD thesis, U.C. Davis, 2003.
\newblock
  \href{http://front.math.ucdavis.edu/math.GT/0307302}{arXiv:math.GT/0307302}.

\bibitem{Casson97}
Andrew~J. Casson.
\newblock The three-sphere recognition algorithm.
\newblock 1997.
\newblock Talk given at MSRI.

\bibitem{CassonGordon87}
Andrew~J. Casson and Cameron~McA. Gordon.
\newblock Reducing {H}eegaard splittings.
\newblock {\em Topology Appl.}, 27(3):275--283, 1987.

\bibitem{ChangDonald91}
Bruce~Randall Donald and David~Renpan Chang.
\newblock On the complexity of computing the homology type of a triangulation.
\newblock In {\em 32nd Annual Symposium on Foundations of Computer Science (San
  Juan, PR, 1991)}, pages 650--661. IEEE Comput. Soc. Press, Los Alamitos, CA,
  1991.
\newblock \href{http://citeseer.ist.psu.edu/chang94computing.html}{URL}.

\bibitem{Gabai87}
David Gabai.
\newblock Foliations and the topology of $3$-manifolds. {I}{I}{I}.
\newblock {\em J. Differential Geom.}, 26(3):479--536, 1987.

\bibitem{GareyJohnson79}
Michael~R. Garey and David~S. Johnson.
\newblock {\em Computers and intractability}.
\newblock W. H. Freeman and Co., San Francisco, Calif., 1979.
\newblock A guide to the theory of NP-completeness, A Series of Books in the
  Mathematical Sciences.

\bibitem{Haken68}
Wolfgang Haken.
\newblock Some results on surfaces in $3$-manifolds.
\newblock In {\em Studies in Modern Topology}, pages 39--98. Math. Assoc. Amer.
  (distributed by Prentice-Hall, Englewood Cliffs, N.J.), 1968.

\bibitem{HassEtAl99}
Joel Hass, Jeffrey~C. Lagarias, and Nicholas Pippenger.
\newblock The computational complexity of knot and link problems.
\newblock {\em J. ACM}, 46(2):185--211, 1999.
\newblock
  \href{http://front.math.ucdavis.edu/math.GT/9807016}{arXiv:math.GT/9807016}.

\bibitem{Hatcher01}
Alan Hatcher.
\newblock Notes on basic 3-manifold topology.
\newblock 2001.
\newblock \href{http://www.math.cornell.edu/~hatcher/3M/3Mdownloads.html}{URL}.

\bibitem{Iliopoulos89}
Costas~S. Iliopoulos.
\newblock Worst-case complexity bounds on algorithms for computing the
  canonical structure of finite abelian groups and the {H}ermite and {S}mith
  normal forms of an integer matrix.
\newblock {\em SIAM J. Comput.}, 18(4):658--669, 1989.

\bibitem{Ivanov01}
S.~V. Ivanov.
\newblock Recognizing the 3-sphere.
\newblock {\em Illinois J. Math.}, 45(4):1073--1117, 2001.
\newblock
  \href{http://www.math.uiuc.edu/~hildebr/ijm/winter01/final/ivanov.html}{URL}.

\bibitem{JacoRubinstein02a}
William Jaco and J.~Hyam Rubinstein.
\newblock 0-efficient triangulations of 3-manifolds.
\newblock 2002.
\newblock
  \href{http://front.math.ucdavis.edu/math.GT/0207158}{arXiv:math.GT/0207158}.

\bibitem{JacoSedgwick03}
William Jaco and Eric Sedgwick.
\newblock Decision problems in the space of {D}ehn fillings.
\newblock {\em Topology}, 42(4):845--906, 2003.
\newblock
  \href{http://front.math.ucdavis.edu/math.GT/9811031}{arXiv:math.GT/9811031}.

\bibitem{King01}
Simon~A. King.
\newblock The size of triangulations supporting a given link.
\newblock {\em Geom. Topol.}, 5:369--398 (electronic), 2001.
\newblock
  \href{http://front.math.ucdavis.edu/math.GT/0007032}{arXiv:math.GT/0007032}.

\bibitem{Rubinstein97}
J.~Hyam Rubinstein.
\newblock Polyhedral minimal surfaces, {H}eegaard splittings and decision
  problems for $3$-dimensional manifolds.
\newblock In {\em Geometric topology (Athens, GA, 1993)}, pages 1--20. Amer.
  Math. Soc., Providence, RI, 1997.

\bibitem{Schleimer01}
Saul Schleimer.
\newblock {\em Almost normal {H}eegaard splittings}.
\newblock PhD thesis, U.C. Berkeley, 2001.
\newblock \href{http://www.math.uic.edu/~saul/Maths/thesis.pdf}{URL}.

\bibitem{Thompson94}
Abigail Thompson.
\newblock Thin position and the recognition problem for ${S}\sp 3$.
\newblock {\em Math. Res. Lett.}, 1(5):613--630, 1994.
\newblock \href{http://www.math.ucdavis.edu/~thompson/thin-paper.aat.pdf}{URL}.

\end{thebibliography}

\end{document}